\documentclass[%
 aip,
 amsmath,amssymb,
 reprint,%
]{revtex4-2}

\usepackage{graphicx}
\usepackage{dcolumn}
\usepackage{bm}

\usepackage[utf8]{inputenc}
\usepackage[T1]{fontenc}
\usepackage{mathptmx}
\usepackage{etoolbox}

\usepackage{lineno}
\pdfoutput=1

\usepackage[pdftex,colorlinks,bookmarksopen,bookmarksnumbered,citecolor=red,urlcolor=red]{hyperref}
\usepackage{amsmath} 	    
\usepackage{amsfonts,amssymb}    	
\usepackage{wrapfig}
\usepackage{siunitx}

\usepackage{mathtools}
\usepackage{graphicx,overpic} 

\newcommand{\Gh}{{\Gamma_h}}
\newcommand{\unablah}{\nabla_{\Gamma_h}}
\newcommand{\rd}{\mathrm{d}}
\newcommand{\bn}{\mathbf{n}}
\newcommand{\R}{\mathbb R}
\newcommand{\bA}{\mathbf A}

\newcommand{\bI}{\mathbf I}

\newcommand{\bP}{\mathbf P}

\newcommand{\bV}{\mathbf V}

\newcommand{\bb}{\mathbf b}

\newcommand{\be}{\mathbf e}
\newcommand{\bj}{\mathbf j}
\newcommand{\bp}{\mathbf p}

\newcommand{\bu}{\mathbf u}

\newcommand{\bv}{\mathbf v}
\newcommand{\bw}{\mathbf w}

\newcommand{\bx}{\mathbf x}

\newcommand{\bbf}{\mathbf f}

\newcommand{\divG}{{\mathop{\,\rm div}}_{\Gamma}}
\newcommand{\gradG}{\nabla_{\Gamma}}
\newcommand{\nablaG}{\nabla_{\Gamma}}
\newcommand{\laplG}{\Delta_{\Gamma}}

\newcommand{\OGamma}{\Omega^\Gamma_h}

\newcommand{\tr}{{\rm tr}}
\DeclareGraphicsExtensions{.pdf,.eps,.ps,.eps.gz,.ps.gz,.eps.Y}

\def\cl {\nonumber \\}
\def\el {\nonumber }

\newtheorem{theorem}{Theorem}

\newcommand{\revo}[1][black]{\textcolor{#1}}
\newcommand{\revt}[1][black]{\textcolor{#1}}

\makeatletter
\def\@email#1#2{%
 \endgroup
 \patchcmd{\titleblock@produce}
  {\frontmatter@RRAPformat}
  {\frontmatter@RRAPformat{\produce@RRAP{*#1\href{mailto:#2}{#2}}}\frontmatter@RRAPformat}
  {}{}
}%
\makeatother
\begin{document}

\preprint{AIP/123-QED}

\title[preprint published as Physics of Fluids 37, 071304 (2025)]{Phase-separated lipid vesicles:
continuum modeling, simulation, and validation}
\author{M. Olshanskii}
\author{A. Quaini}%
 \email{\{maolshan,aquaini\}@central.uh.edu}
\affiliation{ 
Department of Mathematics, University of Houston, 3551 Cullen Blvd, Houston TX 77204, USA
}

\date{ July 18, 2025}

\begin{abstract}
The paper presents a complete research cycle comprising continuum-based modeling, computational framework development, and validation setup to predict phase separation and surface hydrodynamics in lipid bilayer membranes. We \revt{start} with an overview of the key physical characteristics of lipid bilayers, including their composition, mechanical properties, and thermodynamics, and then discuss continuum models of multi-component bilayers. The most complex model is a Navier--Stokes--Cahn--Hilliard (NSCH) type system, describing the 
coupling of incompressible surface fluid dynamics with phase-field dynamics on arbitrarily curved geometries. \revo{This model} is discretized using trace finite element methods, which offer geometric flexibility and stability in representing surface partial differential equations (PDEs). Numerical studies are conducted to examine physical features such as coarsening rates and interfacial dynamics. The computational results obtained from the NSCH model are compared against experimental data for membrane compositions with distinct phase behaviors, demonstrating that including both phase-field models and surface hydrodynamics is essential to accurately reproduce domain evolution observed in epi-fluorescence microscopy. Lastly, we extend the model to incorporate external forces that enable the simulation of vesicles containing cationic lipids, used to enhance membrane fusion. 
\end{abstract}

\maketitle

\section{Introduction}\label{s:intro}

\begin{figure}
 \includegraphics[width=0.3\textwidth]{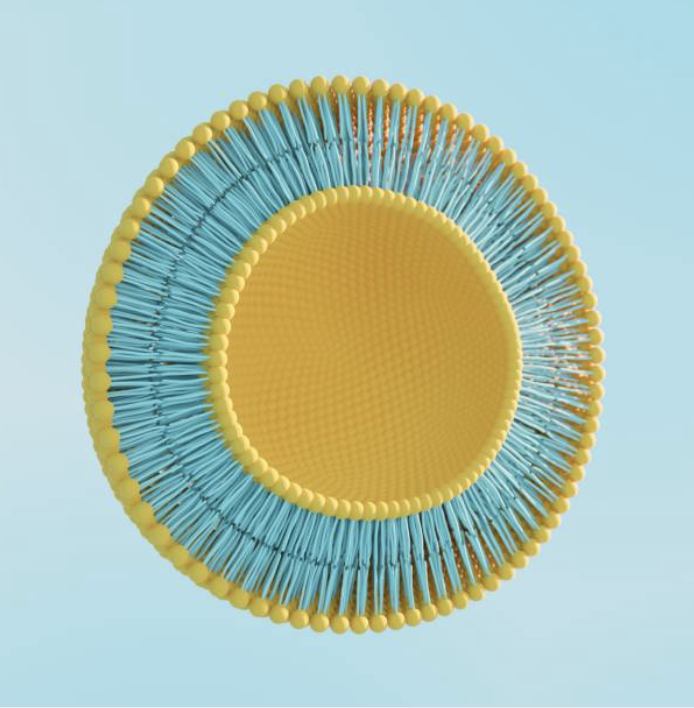}
\caption{Two layers of phospholipids with hydrophilic “heads” in yellow and hydrophobic “tails” in cyan.}\label{fig:lipid_bi}
\end{figure}

Biological membranes separate the interior of a cell from the external environment or create internal boundaries between different cell compartments.
These membranes consist of a bilayer of phospholipids that have a hydrophilic head and a hydrophobic tail.   
In an aqueous environment, membranes form vesicles, 
i.e., bag-like structures that \revt{shield the hydrophobic core from water, minimizing unfavorable energetic interactions}. See Fig.~\ref{fig:lipid_bi}. 
The self-assembly into a stable bilayer structure is governed by hydrophobic interactions, van der Waals forces, and hydrogen bonding~\cite{singer1972fluid,israelachvili1980physical}.
In addition to
lipids, biological membranes contain a mixture of materials including proteins and cholesterol that may decompose into coexisting domains of different phases.

Over the past 20 years, there has been a growing interest in studying phase separation and 
coarsening in multicomponent vesicles. The main reason is its crucial role in a variety of cellular processes\cite{alberts2015essential}. In fact, separation of immiscible liquid phases is likely a factor in the formation of rafts in cell membranes and rafts have been associated with important
biological processes, such as adhesion, signaling, and protein transport. Another reason is that lipid vesicles
provide highly versatile vehicles for intracellular drug delivery and phase separation can enhance the delivery performance.
The formation of artifical lipid bilayers is also central to engineering of model membranes, such as liposomes, which serve as simplified systems for studying membrane behavior~\cite{bangham1965diffusion,szostak2001synthesizing,samimi2019lipid,mazur2017liposomes}.

Multicomponent lipid membranes are liquid-like, resist bending, and are inextensible. Experiments have shown
\revo{that} they exhibit a rich variety of behaviors: spinodal decomposition into distinct surface domains  (liquid-ordered or liquid-disordered), domain coarsening, viscous fingering, vesicle budding, fission, and fusion. 
However, experimental investigations \revt{are} challenging due to the frail nature of giant vesicles, where by giant it is meant that the vesicle diameter is of the order of 10 micron, a size needed for visualization with microscopy. Computational studies can complement experimental studies by allowing \revt{one to observe} dynamics
and gain insights that may not be obtained experimentally.

Monte Carlo methods, dissipative particle dynamics, and molecular dynamics have been used to simulate the dynamics of phase separation and domain formation, vesicle fission and fusion~\cite{marrink2003molecular,laradji2004dynamics,Ehrig_2011,EHRIG201180,venturoli2006mesoscopic,bagatolli2009phase,marrink2019computational}.
However, despite recent advances, these methods still feature high computational costs which limit length and time scales. Continuum methods provide a good modeling alternative to reach larger length and time scales. This paper serves as a widely accessible description of continuum models and efficient numerical methods to approximate the solution. Moreover, it discusses validation against in vitro measurements. 

The remainder of the paper is organized as follows. In Section~\ref{s1}, we provide a brief overview of lipid bilayer composition, mechanical properties, and relevant thermodynamic considerations. Section~\ref{s2} introduces a continuum modeling framework for multi-component lipid bilayers, including surface models for fluid flow, lateral phase separation, and their coupling, as well as elastic deformations and external forces. Section~\ref{s3} presents a finite element based computational method used to simulate the phase-separated vesicle dynamics. In Section~\ref{sec:Exp}, we validate the models against available experimental data, including simulations for the surface fluid–phase separation model and the for system under external forcing. 
Conclusions are drawn in Sec.~\ref{sec:concl}
\revo{and future perspectives are discussed in Sec.~\ref{sec:fp}}.

 \section{Lipid bilayer composition, mechanics, and thermodynamics}\label{s1}

Lipid bilayers of vesicles exhibit a range of structural organizations, including unilamellar and multilamellar arrangements. Unilamellar vesicles consist of a single lipid bilayer enclosing an aqueous core and are widely used in drug delivery due to their ability to encapsulate both hydrophilic and hydrophobic molecules~\cite{lasic1993liposomes,torchilin2005recent,sercombe2015advances,bozzuto2015liposomes}. These vesicles vary in size, ranging from small unilamellar vesicles (SUVs) of 20–100 nm to giant unilamellar vesicles (GUVs) exceeding 10 $\mu$m in diameter. In contrast, multilamellar vesicles contain multiple concentric bilayers, resembling an onion-like structure, and are commonly used in sustained-release drug formulations~\cite{betageri1992drug,gregoriadis1995engineering}.

Under most physiological conditions, lipids diffuse freely in the lateral direction within the bilayer, forming a fluidic membrane embedded in an aqueous environment~\cite{singer1972fluid}. These membranes can exist in different phases, including the liquid-disordered ($L_d$) phase, the liquid-ordered ($L_o$) phase, and the gel ($G$) phase~\cite{cullis1979lipid,gennis2013biomembranes}. The $L_d$ phase is characterized by high lipid mobility and fluidity, typical of unsaturated phospholipids at physiological temperatures. In contrast, the $L_o$ phase, often observed in lipid rafts, exhibits reduced mobility due to the presence of cholesterol and saturated lipids, creating microdomains involved in cell signaling and membrane trafficking. The gel phase, which occurs at lower temperatures, is a rigid and highly ordered state with limited lipid mobility.
In a landmark paper \cite{simons1997functional}, the lipid raft hypothesis was proposed, emphasizing the role of cholesterol and sphingolipids in forming specialized membrane domains. According to this hypothesis, lateral segregation within cell and model membranes arises from lipid phase separations, leading to the formation of biologically relevant entities. Since then, raft formation and lipid bilayer heterogeneity have been subjects of controversy and extensive research, particularly regarding the functional implications of lipid rafts. See, e.g., ~\cite{brown1998functions,ikonen2001roles,munro2003lipid,helms2004lipids,lommerse2004vivo,edidin2003state,binder2003domains,baumgart2007fluorescence,kaiser2009order,sezgin2017mystery,levental2020lipid}.

 Synthetic multicomponent lipid bilayers and vesicles are widely studied as model cell membranes and as components of intracellular delivery carriers.
 Phase-separated lipid bilayers are engineered by carefully selecting lipid compositions and environmental conditions, such as temperature or pH, to promote the coexistence of distinct lipid phases~\cite{veatch2005seeing,veatch2007critical}. The primary approach involves combining lipids with different physicochemical properties—high-melting saturated lipids, low-melting unsaturated lipids, and cholesterol, which modulates lipid packing and phase behavior.
 
 The $L_o$ phase, enriched in saturated lipids and cholesterol, forms more ordered  fluid domains, whereas the $L_d$ phase is dominated by unsaturated lipids. A common strategy to induce phase separation involves ternary mixtures such as sphingomyelin (SM), dioleoylphosphatidylcholine (DOPC), and cholesterol. Sphingomyelin, a saturated lipid, promotes $L_o$ domain formation, while DOPC, an unsaturated lipid, favors the $L_d$ phase. Cholesterol serves as a key regulator, stabilizing the $L_o$ phase by enhancing lipid ordering while maintaining membrane fluidity. Another frequently used lipid combination—dipalmitoylphosphatidylcholine (DPPC), DOPC, and cholesterol—exhibits robust phase separation at physiologically relevant temperatures~\cite{veatch2007critical}.
 
\begin{figure}
\includegraphics[width=0.4\textwidth]{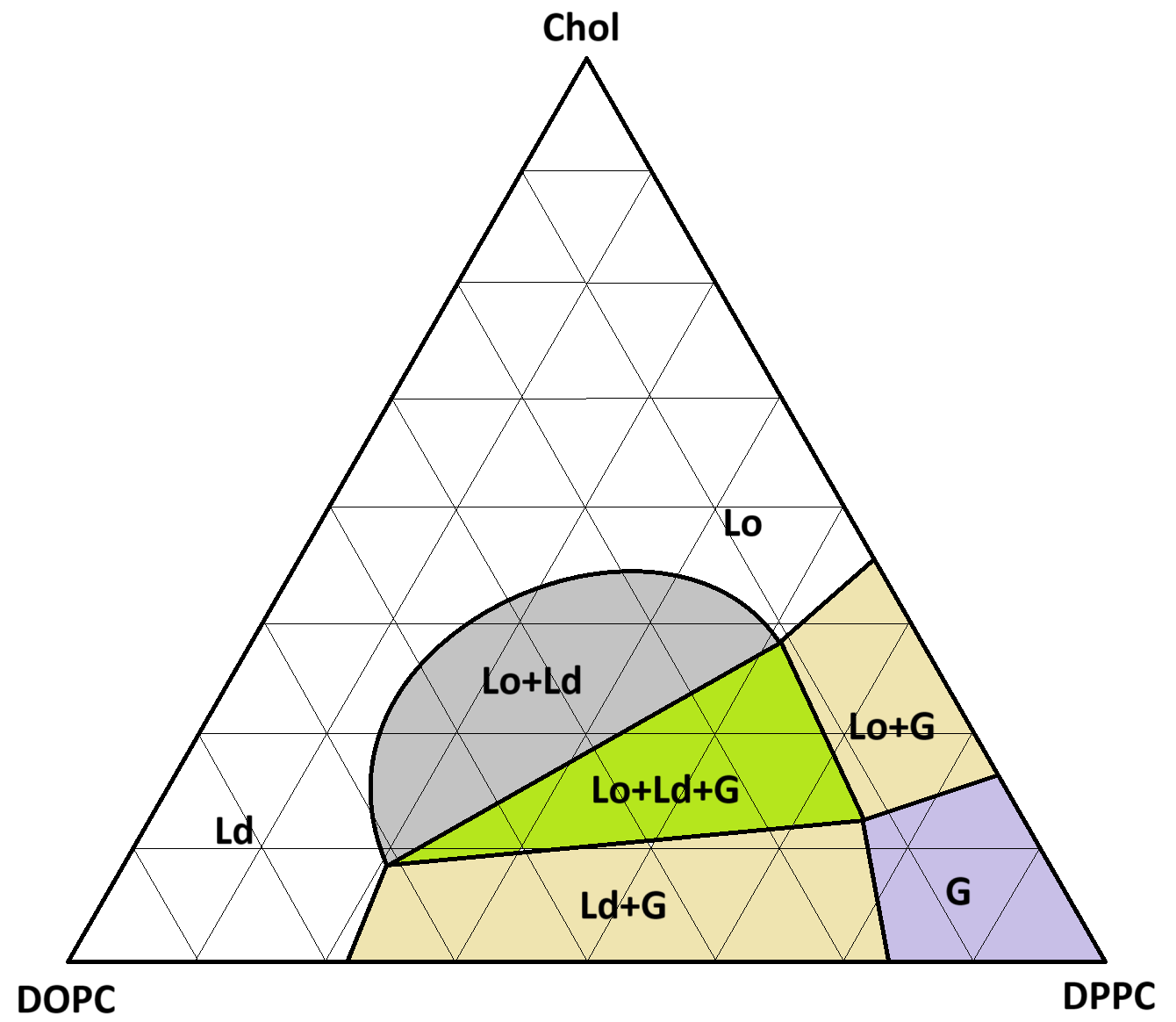}
 	\caption{Schematic phase diagram \revo{qualitatively corresponding to one of the cases analyzed in}~\cite{veatch2007critical}.\label{fig:PD}}
\end{figure}

 Experimental techniques such as fluorescence microscopy are commonly used to construct phase diagrams of lipid bilayers~\cite{veatch2003separation,veatch2007critical,juhasz2009quantitative}.
 A phase diagram is a map 
 of the thermodynamic states of lipid mixtures as a function of variables such as temperature 
 and composition.
 For ternary lipid mixtures, such as DOPC, DPPC, and cholesterol, the phase diagram often displays single-phase regions, two-phase coexistence regions, and three-phase region in which $L_d$, $L_o$, and $G$ phases 
 coexist (see Fig.~\ref{fig:PD}). In the triangular plot, each axis represents the mole fraction of one lipid component. Increasing the cholesterol concentration shifts the system toward the $L_o$ phase, while reducing cholesterol favors $L_d$ behavior. At lower temperatures, some mixtures may also exhibit the gel phase.

 Given the membrane composition and temperature, one can determine the equilibrium composition of each phase and the relative amounts of each phase from the phase diagram. This can be achieved by identifying the corresponding tie-line and applying the lever rule~\cite{chiang2005new,davis2009phase}.
 

The dense packing of lipid molecules and the high energy cost associated with stretching the bilayer make the lipid membrane effectively inextensible~\cite{evans1987physical,rawicz2000effect}. While the fluid membrane is  incompressible  in the lateral direction, it is flexible to deform in the radial direction. The associated bending rigidity of the bilayer depends on its composition and can be measured using a variety of experimental techniques, including fluctuation analysis \cite{helfrich1973elastic,meleard1997bending}, micropipette aspiration \cite{evans1987physical,evans1990entropy}, atomic force microscopy~\cite{dufrene2000advances,garcia2010nanomechanics}, and X-ray or neutron scattering methods~\cite{nagle2000structure}.   

In the context of SUVs used as
drug delivery vehicles, liposomes that contain cationic lipids \cite{hui1996role,almofti2003cationic}, 
such as 1,2-dioleoyl-3-trimethylammonium-propane (DOTAP), are known 
for their high fusogenicity \cite{simberg2004dotap}, i.e., the enhanced ability to fuse with the target cell and 
deliver the encapsulated cargo directly into the cytoplasm.
However, while cationic lipids are typically non-toxic at lower concentrations, 
concerns arise regarding their toxicity when used at higher concentrations \cite{lv2006toxicity}. Hence, 
membrane phase separation has been considered to design
delivery liposomes that offer both 
high fusogenicity and low toxicity. For example, 
phase-separation 
in mixture containing DOTAP, DOPC, DPPC, and 
cholesterol was studied 
to modulate surface density of DOTAP on liposomes and hence their fusogenicity \cite{wang2024}. 

In summary, lipid bilayers, also known as fluid or lipid membranes, are self-organized, amphipathic lipid sheets that are two molecules thick and immersed in an aqueous environment. They form the membranes of all living cells, while engineered bilayers are widely used in biomedical applications. The key physical properties essential for both their function and mathematical modeling include lateral fluidity, out-of-plane elasticity, and phase transitions, which lead to the coexistence of different phases depending on membrane composition and external conditions. The thermodynamic and mechanical characteristics of lipid bilayers, such as  density, viscosity, bending moduli, and phase diagrams, have been quantitatively assessed through various \textit{in vitro} measurements and extensively documented in the literature.

\section{Continuum models of multi-component bilayers} \label{s2}

To model the key physical properties of lipid vesicles, we start with a fundamental observation: lipid bilayers are extremely thin (on the order of a few nanometers) relative to their lateral dimensions.
In fact, the diameter of a mammalian cell is about 10 micrometers and the typical diameter of a liposome is about 100 nanometers. Thus, it is reasonable to approximate the lipid vesicle with a closed, sufficiently smooth surface $\Gamma$ embedded in $\mathbb{R}^3$. Since the continuum models 
will be posed on $\Gamma$, let us introduce some basic elements of tangential calculus on embedded surfaces.

The outward pointing unit normal on $\Gamma$ is denoted by $\bn$, and the orthogonal projection on the tangential plane is given by $\bP=\bP(\bx):= \bI - \bn(\bx)\bn(\bx)^T$, $\bx \in \Gamma$. In a neighborhood $\mathcal{O}(\Gamma)$  of $\Gamma$ the closest point projection $\bp:\,\mathcal{O}(\Gamma)\to \Gamma$ is well defined.  
For a scalar function $p:\, \Gamma \to \mathbb{R}$ or a vector function $\bu:\, \Gamma \to \mathbb{R}^3$  we define their extensions from $\Gamma$ to its neighborhood $\mathcal{O}(\Gamma)$ along the normal directions as
$$
p^e=p\circ \bp\,:\,\mathcal{O}(\Gamma)\to\mathbb{R}, \quad
\bu^e=\bu\circ \bp\,:\,\mathcal{O}(\Gamma)\to\mathbb{R}^3.
$$
On $\Gamma$, it holds $\nabla p^e= \bP\nabla p^e$  and  $\nabla\bu^e=\nabla\bu^e\bP$, with $\nabla \bu:= (\nabla u_1~ \nabla u_2 ~\nabla u_3)^T \in \mathbb{R}^{3 \times 3}$ for vector functions $\bu$. Hence, the surface gradient and covariant derivatives on $\Gamma$ are defined as $\nablaG p=\bP\nabla p^e$ and  $\nabla_\Gamma \bu:= \bP \nabla \bu^e \bP$. Note that the definitions  of surface gradient and covariant derivatives are  independent of a particular smooth extension of $p$ and $\bu$ off $\Gamma$.
We consider normal extensions because they are  
convenient for the error analysis of the computational method (see Sec.~\ref{s3}).

If vector function $\bu:\, \Gamma \to \mathbb{R}^3$ represents a velocity field, 
we define the surface rate-of-strain tensor \cite{GurtinMurdoch75} on $\Gamma$ as
\begin{equation} \label{strain}
 E_s(\bu):= \frac12 \bP (\nabla \bu +\nabla \bu^T)\bP = \frac12(\nabla_\Gamma \bu + \nabla_\Gamma \bu^T).
 \end{equation}
The surface divergence operator of vector $\bu$ is given by 
\[
 \divG \bu  := \tr (\gradG \bu).
               \]
We also define the surface divergence operator of
a tensor $\bA: \Gamma \to \mathbb{R}^{3\times 3}$ as:
\[
 \divG \bA  := \left( \divG (\be_1^T \bA),\,
               \divG (\be_2^T \bA),\,
               \divG (\be_3^T \bA)\right)^T,
               \]
with $\be_i$ the $i$th basis vector in $\R^3$.

Finally, for a smooth field $f: \mathbb{R}^n \rightarrow \mathbb{R}$, we denote by $\dot{f}$
the material derivative, i.e., the derivative along
material trajectories in a velocity field
$\bu$:
\begin{equation*}
    \dot{f} = \frac{\partial f}{\partial t} + (\nabla^T f) \cdot \bu.
\end{equation*}
Derivative $\dot{f}$ is a tangential derivative for $\Gamma$ and so it depends only
on the surface values of f on $\Gamma$.
For a vector field $\bv$, $\dot{\bv}$ is defined component-wise.

\subsection{Surface model for fluid flow}\label{sec:NS}

Let $\rho$ be a density distribution on $\Gamma$, \revo{i.e., $\rho$ is the surface density of the membrane}, and let 
$\bu:\, \Gamma \to \R^3$ denote the smooth velocity field
of the density flow on $\Gamma$.
Let us write $\bu$ as the sum of its tangential and normal components:
\begin{equation}\label{u_T_N}
\bu = \overline{\bu}+(\bu\cdot\bn)\bn =  \overline{\bu} + u_N\bn,\quad \overline{\bu}\cdot\bn=0.
\end{equation}
Since we consider $\Gamma$ at shape equilibrium, we have 
$\bu\cdot\bn = u_N = 0$. 

To state inextensibility and conservation of mass, let us introduce a material subdomain $\gamma \subset \Gamma$ and Leibniz rule, which written for a smooth function $f: \Gamma \to \R$ reads:
\begin{equation}\label{eq:Leibniz}
    \frac{\partial}{\partial t} \int_\gamma f\,ds = \int_\gamma (\dot{f} + f \divG \bu)\,ds.
\end{equation}

Inextensibility means that
\begin{equation*}
    0 = \frac{\partial}{\partial t} \int_\gamma 1\,ds = \int_\gamma \divG \bu\,ds,
\end{equation*}
where we have applied \eqref{eq:Leibniz}. Since $\gamma$ is arbitrary, we obtain 
\begin{equation}\label{strongform-2}
    \divG \bu =0 \quad \text{on}~\Gamma.
\end{equation}

From conservation of mass, by applying \eqref{eq:Leibniz} and \eqref{strongform-2}, we obtain:
\begin{equation*}
    0 = \frac{\partial}{\partial t} \int_\gamma \rho\,ds = \int_\gamma \dot{\rho}\,ds,
\end{equation*}
which means $\dot{\rho} = 0$ on $\Gamma$ since $\gamma$ is arbitrary. Thus, if $\rho |_{t = 0} =$ const, then $\rho$ is constant for all $t > 0$. 

The conservation of momentum for $\gamma$ reads:
\begin{equation}\label{eq:mom}
    \frac{\partial}{\partial t} \int_\gamma \rho \bu \,ds = \int_{\partial \gamma} \boldsymbol{\sigma}_\Gamma \boldsymbol{\nu} \,dl + \int_\gamma \bb \,ds,
\end{equation}
where $\boldsymbol{\nu}$ is the normal to
$\partial \gamma$, 
$\bb$ denotes the area forces on $\gamma$,
and $\boldsymbol{\sigma}_\Gamma$ is the Boussinesq–Scriven surface stress tensor, which is given by:
\begin{equation*}
    \boldsymbol{\sigma}_\Gamma = -p \bP + 2 \eta E_s(\bu),
\end{equation*}
once inextensibility is taken into account. Here, 
$\eta$ is the dynamic viscosity of the fluid, \revo{i.e., the membrane shear viscosity}.
By applying Leibniz rule to the left-hand side of \eqref{eq:mom} and accounting for \eqref{strongform-2} and conservation of mass (i.e.,  $\dot{\rho} = 0$ on $\Gamma$), we get:
\begin{align}
    &\frac{\partial}{\partial t} \int_\gamma \rho \bu \,ds =  \int_\gamma (\dot{\rho} \bu + \rho \dot{\bu} + \rho \bu \divG \bu) \,ds
    = \int_\gamma \rho \dot{\bu} \,ds \cl
    & \quad = \int_\gamma \rho \left( \frac{\partial \bu}{\partial t} + (\gradG \bu) \bu \right) \,ds. \el
\end{align}
By applying the divergence theorem to the right-hand side of \eqref{eq:mom}, we obtain:
\begin{equation*}
    \int_{\partial \gamma} \boldsymbol{\sigma}_\Gamma \boldsymbol{\nu} \,dl + \int_\gamma \bb \,ds = \int_\gamma  (\divG(- p \bP + 2 \eta E_s(\bu)) + \bb)\,ds.
\end{equation*}
So, given that $\gamma$ is arbitrary and $\Gamma$ is stationary, from \eqref{eq:mom} we get:
\begin{align} 
  \rho  \frac{\partial \bu}{\partial t} + \rho (\gradG \bu) \bu - \bP \divG (2 \eta E_s(\bu)) +\nabla_\Gamma p =  \bb \quad \text{on}~\Gamma.  \label{strongform-1}
\end{align}

The surface Navier--Stokes problem reads: 
Find a vector field $\bu:\, \Gamma \to \R^3$, with $\bu\cdot\bn =0$, such that eq.~\eqref{strongform-2},\eqref{strongform-1} 
hold. For this problem, the pressure field is defined up a hydrostatic constant. Further, if the flow is stationary (i.e., \revt{when} $\dot{\bu}$ is negligible), all tangentially rigid surface fluid motions (i.e., with  $E_s(\bu)=0$), also known as Killing fields \cite{sakai1996riemannian},
are in the kernel of the differential operators at the left-hand side of \eqref{strongform-1}.
Thus, for consistency the following condition is required:
 \begin{equation}\label{constr}
 \int_\Gamma \bb \cdot \bv\,ds=0
 \end{equation}
for all smooth tangential vector fields $\bv$ s.t.~$E_s(\bv)=\mathbf{0}$.
This condition  is necessary for the well-posedness
of problem \eqref{strongform-2},\eqref{strongform-1} when the flow is stationary.

The reader interested in 
the derivation of the Navier--Stokes equations  for evolving fluidic interfaces is referred to, e.g., \cite{Jankuhn1}.

\subsection{Surface model for phase separation}\label{sec:CH}

On $\Gamma$, we consider a heterogeneous mixture of two species with
surface fractions $c_i = S_i/S$, $i = 1, 2$, where $S_i$
are the surface area occupied by the components and $S$ is the surface area of $\Gamma$.
Since $S = S_1 + S_2$, we have $c_1 + c_2 = 1$. Let $c_1$ be the representative surface fraction, i.e.,~$c = c_1$.
Moreover, let $m_i$ be the mass of component $i$ and $m$ is the total mass.
Notice that density of the mixture can be expressed as
$\rho= \frac{m}{S} = \frac{m_1}{S_1} \frac{S_1}{S} + \frac{m_2}{S_2} \frac{S_2}{S}$.
Phase separation in this two component system
can be modeled by the Cahn--Hilliard (CH) equation \cite{Cahn_Hilliard1958,CAHN1961}.

In order to describe the evolution of $c(\bx, t)$, we consider the conservation law:
\begin{align}\label{eq:evol_c}
\rho \frac{\partial c}{\partial t} +  \divG \bj = 0 \quad \text{on}~\Gamma  \times (0, T],
\end{align}
where $\bj$ is a diffusion flux. The flux $\bj$ is defined according to Fick's law:
\begin{align}\label{eq:flux}
\bj = - M \gradG \mu \quad \text{on}~\Gamma, \quad \mu = \frac{\delta f}{\delta c},
\end{align}
where $M$ is the so-called mobility coefficient (see \cite{Landau_Lifshitz_1958}) and
$\mu$ is the chemical potential, which is defined as the functional derivative of the total
specific free energy $f$ with respect to the concentration $c$. Many analytic studies, as well as numerical simulations, assume mobility is constant. Another popular choice for numerical studies is 
\begin{equation}\label{eq:mob}
 M = M(c) = D\,c(1 - c),   
\end{equation}
 which \revt{is} called degenerate, since it is not strictly positive. Only a few authors consider more complex
mobility functions; see, e.g., \cite{PhysRevE.60.3564}. We note that there is virtually no study on the
appropriate mobility function for lateral phase separation in biological membranes. \revo{See Sec.~\ref{sec:est} for details on the setting of mobility
for our validation studies reported in Sec.~\ref{sec:NSCH_val} and \ref{sec:charge_val}.}

\revt{We} now introduce
the total specific free energy:
\begin{align}\label{eq:total_free_e}
f(c) = \frac{1}{\epsilon}f_0(c) + \frac{1}{2} \epsilon | \gradG c |^2.
\end{align}
where $f_0(c)$ is the free energy per unit surface, which 
must be a non-convex function of $c$ for phase separation to occur. A common choice for $f_0$ is \revo{the Ginzburg--Landau double-well potential:}
\begin{align}\label{eq:f0}
f_0(c) = \frac{1}{4} c^2(1 - c)^2.
\end{align}
The second term in \eqref{eq:total_free_e} represents
the interfacial free energy based on the surface fraction gradient. The physical meaning of parameter $\epsilon$ in 
\eqref{eq:total_free_e} is the thickness of a layer where
thermodynamically
unstable mixtures are stabilized by a gradient term in the energy.

By combining eq.~\eqref{eq:evol_c}, \eqref{eq:flux}, and \eqref{eq:total_free_e}, we obtain
the surface CH equation:
\begin{align}\label{eq:CH}
\rho \frac{\partial c}{\partial t} -  \divG \left(M \gradG \left(\frac{1}{\epsilon} f_0' - \epsilon \laplG c\right)\right) = 0 \quad \text{on}~\Gamma.
\end{align}
With the specified choices for mobility and
specific free energy, eq.~\eqref{eq:CH} is a nonlinear fourth-order equation.
Since one might be interested in avoiding higher order 
spatial derivatives in view of a numerical algorithm to solve the problem, it is common to rewrite eq.~\eqref{eq:CH} in mixed form:
\begin{align}
&\rho \frac{\partial c}{\partial t} -  \divG \left(M \gradG \mu \right)  = 0 \quad \text{on}~\Gamma \label{eq:sys_CH1}, \\
&\mu = \frac{1}{\epsilon} f_0' - \epsilon \laplG c \quad \text{on}~\Gamma. \label{eq:sys_CH2}
\end{align}
Problem \eqref{eq:sys_CH1}--\eqref{eq:sys_CH2} is a coupled system of nonlinear PDEs posed on $\Gamma$.

We note that the CH model is conservative, as it comes from conservation law \eqref{eq:evol_c}.
Experimentally, a number of molecules are known to preferentially partition into one of lipid phases on phase-separated vesicles.
Examples include membrane proteins caveolin-3, peripheral myelin protein 22 \cite{Baumgart3165,Schlebach2016}
and membrane dyes \cite{Baumgart3165,BAUMGART2007}. In these cases, the CH problem provides the correct model. Similarly, the conservative model seems to be suited to describe membrane separation in bacteria as well as mammalian cells. 
The use of a non-conservative model (e.g., the Allen--Cahn equation \cite{ayton2005coupling,elliott2010modeling,Wang2008,ElliottStinner,elliott2013computation}) may be justified
when phase separation induces ``high'' curvature to the membrane that leads to vesicle budding
or, basically, formation and separation of individual vesicles from the original parent vesicle
\cite{Baumgart_et_al2003,HURLEY2010875}.

The reader interested in the surface
Cahn--Hilliard problem posed on evolving surfaces is referred to, e.g., \cite{Jankuhn1}.

\subsection{Coupled surface model for fluid flow and phase separation}\label{sec:NSCH}

In this section, we combine the models presented in Sec.~\ref{sec:NS} and \ref{sec:CH}, using the notation introduced therein.
The classical phase-field model to describe the flow of two immiscible, incompressible, and Newtonian fluids
is called \emph{Model H} \cite{RevModPhys.49.435}, which is a Navier–Stokes–Cahn–Hilliard (NSCH) system. One of the fundamental assumptions for Model H
is that the densities of both components are matching. Several extensions have been proposed
to account for the case of non-matching densities. See, e.g., \cite{Abels2012,Aki2014,BOYER200241,DING20072078,DONG20125788,GONG201720,Gong2018,Lowengrub1998,shokrpour2018diffuse}.
For example, the thermodynamically consistent generalization of Model H in \cite{Abels2012} allows for non-matching densities $\rho_1, \rho_2$ for the two phase but forces the density of the mixture to obey:
\begin{align}
\rho= \rho(c) = \rho_1 c+ \rho_2 (1-c), \label{rho_c}
\end{align}
In \cite{palzhanov2021decoupled}, we considered
the model from \cite{Abels2012}
and extended it more general relations than \eqref{rho_c}.
This was motivated by the fact that, depending on the choice of $f_0(c)$ and because of numerical errors, $c$ may not be constrained in $[0,1]$ and so $\rho$ in \eqref{rho_c} may take physically meaningless (i.e., negative) values. Below, we present the model \revt{proposed in} \cite{palzhanov2021decoupled}.

Let $\rho_1\ge\rho_2$ and assume that $\rho$ is a smooth monotonic function of $c$, i.e., $\frac{d\rho}{dc}\ge 0$ (for $\rho_1\ge\rho_2$), and so we can set
\begin{equation}\label{monot}
\frac{d\rho}{dc} = \theta^2.
\end{equation}
Let $\sigma_\gamma$ be line tension and 
$\bbf$ a given force vector. 
The surface NSCH model from \cite{palzhanov2021decoupled} reads:
\begin{align} 
  &\underbrace{\rho \frac{\partial \bu}{\partial t}  + \rho(\nabla_\Gamma\bu)\bu}_{\revo{\text{inertia}}} - \underbrace{\bP\divG(2\eta E_s(\bu)) +\nabla_\Gamma p}_{\revo{\scriptsize\begin{array}{c}\text{viscous and in-plane pressure}\\ \text{gradient forces}\end{array}}}\hspace{10ex}\cl
  & \quad \quad \quad=  \underbrace{-\sigma_\gamma c  \nablaG \mu}_{\revo{\text{line tension}}} + \underbrace{M \theta(\nabla_\Gamma(\theta\bu)\,)\gradG \mu}_{\revo{\text{chemical momentum flux}}} + \bbf,  \label{grache-1m} \\
& \underbrace{\divG \bu  =0,}_{\revo{\text{membrane inextensibility}}}  \label{gracke-2}\\
&\underbrace{\frac{\partial c}{\partial t}  +\divG(c\bu)}_{\revo{\text{transport of phases}}}-  \underbrace{\divG \left(M \gradG \mu \right)}_{\revo{\scriptsize\begin{array}{c}\text{phase masses exchange}\\ \text{Fick's law}\end{array}}}  = 0  \label{gracke-3}, \\
&\mu = \underbrace{\frac{1}{\epsilon} f_0' - \epsilon \laplG c,}_{\revo{\scriptsize\begin{array}{c}\text{mixture free}\\ \text{energy variation}\end{array}}} \label{gracke-4}
\end{align}
Note that since $\Gamma$ is stationary, system \eqref{grache-1m}--\eqref{gracke-4}
is fully tangential. 
The only term \revt{by} which model \eqref{grache-1m}-\eqref{gracke-4} differs from the model in \cite{Abels2012} is the middle term at the right-hand side in eq.~\eqref{grache-1m}. This term, crucial for thermodynamic consistency, can be interpreted as  additional momentum flux due to diffusion of
the components driven by the gradient of the chemical potential. Note that when $\rho_1 = \rho_2$, i.e., $\rho$ does not depend on $c$, this term vanishes.
The other terms in \eqref{grache-1m}-\eqref{gracke-4} that are not present in either \eqref{strongform-2}--\eqref{strongform-1} or \eqref{eq:sys_CH1}--\eqref{eq:sys_CH2}
are the first term at the right-hand side in eq.~\eqref{grache-1m}, which models line tension, and the second term at the left-hand
side in eq.~\eqref{gracke-3}, which contributes to the transport of phases.

We note that temperature does not appear in eq.~\eqref{grache-1m}--\eqref{gracke-4},
which describe the evolution of phases and coupled surface flow independently of what initiates phase separation. As such, 
the same model could be used if phase separation is triggered by, e.g., pH \cite{Bandekar2012} or temperature.

Other thermodynamic consistent extensions of  the model in \cite{palzhanov2021decoupled} for a generic smooth $\rho(c)$, with no monotonicity assumption, were analyzed in \cite{abels2016weak,abels2019existence} in terms of well-posedness analysis. Those extensions introduce more terms in the momentum equation, so they are slightly more cumbersome for computations and numerical analysis.

\subsection{Surface models for elasticity}
Following the success of the Canham--Helfrich energetic approach~\cite{canham1970minimum,helfrich1973elastic} in describing the equilibrium shapes of red blood cells, standard elasticity models for lipid membranes are based on the principle of bending energy minimization. The bending energy of the membrane is commonly defined by the Willmore functional:
\begin{equation} \label{Willmore}
	H = \frac{c_\kappa}{2}\int_{\Gamma}(\kappa - \kappa_0)^2\,ds + c_K\int_{\Gamma}K\,ds,
\end{equation}
or one of its variants~\cite{seifert1991shape,seifert1997configurations}.
Here, $K$ is the Gaussian curvature, and the material parameters $c_\kappa > 0$, $c_K > 0$, and $\kappa_0$ represent the bending rigidity, Gaussian bending rigidity, and spontaneous curvature, respectively.

To derive the elastic forces arising from variations in the bending energy, one may apply the principle of virtual work, yielding
\[
\int_{\Gamma(t)} \bb^{\rm elst} \cdot \bv \, ds = -\left.\frac{dH}{d\Gamma}\right|_{\bv},
\]
where $\left.\frac{dH}{d\Gamma}\right|_{\bv}$ denotes the variation of the energy functional under an (infinitesimal) displacement of $\Gamma$ defined by the vector field $\bv$.
The shape derivative of $H$ can be computed~\cite{willmore1996riemannian} and takes the form:
\begin{equation}\label{BE1}
\left.\frac{dH}{d\Gamma}\right|_{\bv} = c_\kappa \int_{\Gamma(t)} (-\Delta_\Gamma \kappa - \tfrac{1}{2} \kappa^3 + 2K\kappa)(\bv \cdot \bn) \, ds.
\end{equation}
It follows from \eqref{BE1} that the release of bending energy generates a force in the normal direction to the surface:
\begin{equation} \label{bN}
	\bb^{\rm elst} = c_\kappa \left( \Delta_\Gamma \kappa + \tfrac{1}{2} \kappa^3 - 2K\kappa \right) \bn.
\end{equation}

\revo{Another normal force acting on membrane is the surface tension. The corresponding term would appear in the governing equations if one allows  radial motions, i.e., when $\bu\cdot\bn\neq 0$. See, e.g., \cite{Jankuhn1}.}
Since we consider vesicles of fixed shape, we assume that the bending force is balanced by the surface tension and osmotic pressure forces acting on $\Gamma$. \revo{These forces define the shape of the surface $\Gamma$, so they do not appear in the governing equations \eqref{grache-1m}--\eqref{gracke-4}.}

It is important to account for $\bb^{\rm elst}$ when the deformation dynamics of the membrane is of interest~\cite{arroyo2009relaxation,hu2007continuum,torres2019modelling,reuther2020numerical}, or when determining the equilibrium shape of a vesicle is part of the modeling problem~\cite{deuling1976curvature,jenkins1977static,deserno2015fluid,seifert1991shape,seifert1997configurations,seifert1995morphology,olshanskii2023equilibrium}.

\subsection{Estimation of modeling parameters}\label{sec:est} 

The NSCH model \eqref{grache-1m}--\eqref{gracke-4} assumes the coexistence of lipid phases in thermodynamic equilibrium (i.e., variations of the temperature are 
thermodynamically insignificant after the phase separation is initiated if temperature is the trigger). In this case, existing phase diagrams are applied to determine the phase ratio in the initial condition, since this ration is conserved by the system. 
Key physical properties governing multiphase bilayer evolution include membrane density, viscosity, and interphase tension forces. The  densities $\rho_i$, $i = 1, 2$, can be estimated from molecular weight and molecular surface area for each phase~\cite{wang2022lipid}, while viscosities $\eta_i$, $i = 1, 2$, which cannot be measured with a traditional viscometer, are assessed through bilayer responses to external forces~\cite{sakuma2020viscosity}. 
For instance, membrane viscosity is estimated by measuring translational diffusion coefficients of tracer particles embedded in membranes, including lipids~\cite{oradd2004lateral}, fluorescent probes~\cite{yguerabide1982lateral}, labeled proteins~\cite{gambin2006lateral}, or membrane-linked particles~\cite{hormel2014measuring}. Hydrodynamic models are then used to derive membrane viscosity from these measurements. Additionally, several studies~\cite{heftberger2015situ,kollmitzer2013monolayer,kuzmin2005line} have estimated interphase line tension forces, from which the line tension coefficient $\sigma_\gamma$ can be estimated.

In addition, the system \eqref{grache-1m}--\eqref{gracke-4} depends on two key modeling parameters: $\epsilon$, which defines the width of the transition layer between ordered and disordered phases in the free energy, and the diffusivity coefficient $D$ in the mobility term \eqref{eq:mob}. While both $\epsilon$ and $D$ correspond to thermodynamic properties, their direct measurement is nontrivial. In particular, $D$ governs the rate of change of $c$ based on free energy fluctuations, rather than molecular diffusion via Brownian motion. In \eqref{gracke-4}, the timescale scales linearly with $D$, while $\epsilon$ determines the relative duration of the rapid decomposition phase versus the slower coarsening process.
Given the uncertainty in $D$ and $\epsilon$, a data-driven approach was proposed in \cite{zhiliakov2021experimental}, using backward optimization to calibrate these parameters against observed \textit{in vitro} pattern dynamics. There, $D$ was estimated in the range $10^{-5}$–$2.5 \times 10^{-5}$ cm$^2$/s depending on membrane composition. For $\epsilon$, a value of \SI{0.1}{\micro\metre} provided good agreement with experiments, though this overestimates the $\sim$5 nm transition width found in \cite{Risselada2008}. This discrepancy reflects the resolution limits of the discrete continuum model, effectively broadening the interfacial region where tension forces act while preserving the resulting momentum.

\subsection{External forces}\label{sec:ex_forces} 
In this section, we consider a particular forcing term $\bbf$ in \eqref{grache-1m} that is relevant to the study 
of phase-separation in liposomes (referred to as SUV) that contain cationic lipids. See Sec.~\ref{s1}. 

\begin{figure}[htb!]
	\centering
	\begin{overpic}[width=.50\textwidth,grid=false]{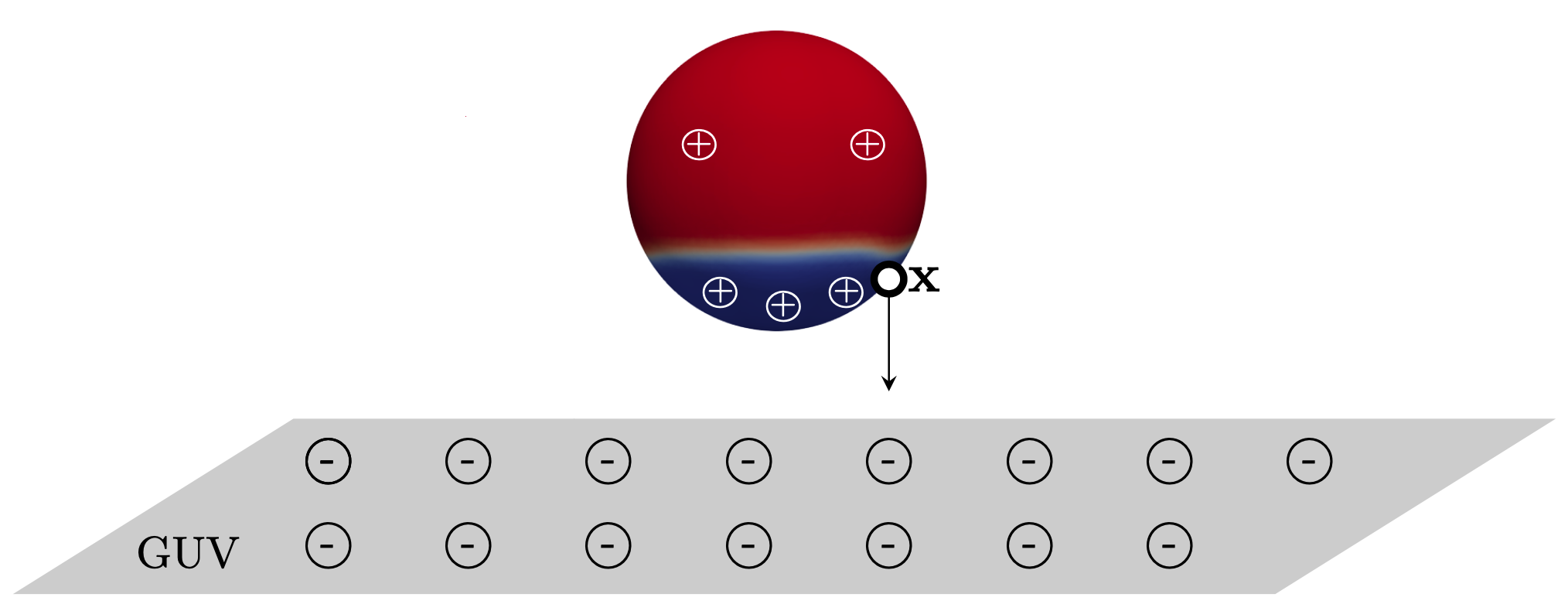}
    \put(25,25){\footnotesize{SUV}}
\end{overpic}
	\caption{Simulation set-up: GUV, represented as a plane, and a 
	phase-separated SUV.  The $L_o$ phase in 
	the phase-separated SUV is colored in red, while the $L_d$ phase is blue. }
	\label{fig:force}
\end{figure}

With the aim of examining the effect of positively charged lipids on the fusogenicity of the SUVs, 
let us assume that a SUV is in the vicinity of a GUV,
which represents a model target cell
and has a slight negative charge \cite{chibowski2016,MAKINO1991175,molecules25122780}.
Both SUV and GUV are immersed in a NaCl solution typically used in lab experiments. In this case,
a relevant forcing term is the electrostatic force per unit surface area acting on the SUV, denoted with $\bbf_e$.
Because the GUVs are significantly larger than the SUVs, 
the curvature of a GUV is negligible {at the scale given by the size} of an SUV. 
See Fig.~\ref{fig:force}.
Hence, we will approximate a GUV with a plane
for the computation of the electrostatic force $\bbf_e$. Therefore, the electric field $\textbf{E}$ generated by a GUV
{can be (locally) computed by}:
\begin{equation}\label{eq:E}
\boldsymbol{E} = \frac{\sigma}{2\varepsilon_0}, 
\end{equation}
where $\sigma$ is the GUV surface charge density and $\varepsilon_0$ is the vacuum permittivity ($8.85\cdot10^{-12}$ \SI{}{\farad\per\meter}). 
Since the value of $\sigma$ cannot be measured, it is estimated from a linear approximation of
Grahame's formula \cite{interfacesbook}, which is valid in low-potential situations:
\begin{equation}
	\sigma   \approx \varepsilon\cdot \varepsilon_0 \cdot\kappa \cdot \Psi_0 , \quad  \Psi_0 = \frac{\zeta}{\exp(-\kappa\cdot x)},
	\label{grahame}
\end{equation}
where $\varepsilon$ is the relative permittivity of water (about 80 at 20\textdegree{}C), 
$\kappa$ is the Debye length parameter for a  NaCl solution (10/7 $\text{nm}^{-1}$),  
$\Psi_0$ is the surface potential \cite{chibowski2016zeta}, 
$x$ is the slip plane (\SI{0.24}{\nano\meter}), 
and $\zeta$ is the zeta potential. 
The zeta potentials for GUVs and SUVs of given compositions can be measured
experimentally. See, e.g., \cite{wang2024}.

Once the electric field $\textbf{E}$ is computed, $\bbf_e$ is given by
$\bbf_e(\bx) = \textbf{E} q(\bx)$, where $q$ is a point charge located at
$\bx$ on an SUV. See Fig.~\ref{fig:force}.
A point charge on an SUV cannot be measured, 
hence we resort to 
an approximation: with the 
surface charge density \eqref{grahame} for an SUV
of given composition obtained from the measured zeta potential, we get the total attraction force density and we distribute it proportionally to the SUV surface. See Sec.~\ref{sec:charge_val} for more details
on this approximation.

\section{Computational methods}\label{s3}
Numerical methods for solving PDEs on surfaces have been an active area of research over the past two decades and the subject of several review papers~\cite{dziuk2013finite,olshanskii2017trace,bonito2020finite}. The study of finite element methods (FEMs) for PDEs on general surfaces dates back to Dziuk's work~\cite{Dziuk88}, which focused on the Laplace–Beltrami equation on a stationary surface \(\Gamma\), approximated by a family of triangulations $\{\Gamma_h\}$ with vertices on $\Gamma$. Initially introduced for linear elements, this approach was later extended to higher-order elements~\cite{Demlow09} and adaptive FEM with \textit{a posteriori} error estimators~\cite{Demlow06}. 

To avoid remeshing and leverage the implicit surface definition via level-set functions, Bertalmio et al.~\cite{bertalmio2001variational} proposed extending the PDE from the surface to a higher-dimensional domain in $\mathbb{R}^3$. This approach allows solving the PDE on a background mesh not aligned with the surface. It has been explored for finite difference methods~\cite{AS03,Greer,xu2006level,XuZh}, including applications to moving surfaces, and has been combined with FEM~\cite{burger2009finite,DziukElliot2010}. 
A related technique is the closest point method~\cite{macdonald2009implicit,ruuth2008simple}, which extends the surface problem to a neighborhood using the closest point operator and discretizes it with Cartesian finite differences. Unlike Dziuk’s method, these approaches fall into the category of unfitted or immersed interface discretizations.  

The TraceFEM~\cite{olshanskii2009finite} discussed in this article is another example of an unfitted or immersed interface method and is closely related to CutFEM~\cite{burman2015cutfem}, which originally emerged as an unfitted FEM for interface problems. These methods eliminate the need for surface-aligned meshes, offering greater flexibility in handling complex geometries. In TraceFEM, the surface is allowed to freely "cut" through the background (ambient) mesh, and finite element function traces from the bulk space are used to approximate functions defined on the surface. Next, we introduce the necessary notation and
definitions to apply TraceFEM for the numerical solution of the models presented in Sec.~\ref{s2}.

Let $\mathcal{T}_h$ be a tetrahedral triangulation of the domain $\Omega \subset \mathbb{R}^3$ that contains $\Gamma$. This triangulation is assumed to be regular, consistent, and stable~\cite{Braess}; it serves as the background mesh for the TraceFEM.  
We assume an approximation $\Gamma_h$ of $\Gamma$ such that the integration of finite element (FE) functions over $\Gamma_h$ is feasible. Possible constructions of $\Gamma_h$ include the zero level set of $P^1$ FE level set functions or higher-order implicit surface reconstructions, as described in~\cite{lehrenfeld2016high}.  

The active set of tetrahedra $\mathcal{T}_h^\Gamma \subset \mathcal{T}_h$ is defined as  
\[
\mathcal{T}_h^\Gamma = \{\, T \in \mathcal{T}_h \,:\, \text{meas}_2(\Gamma_h \cap T) > 0\, \}.
\]  
The domain formed by the tetrahedra in $\mathcal{T}_h^\Gamma$ is denoted by $\omega_h$. In the TraceFEM, only the background degrees of freedom associated with the tetrahedra in $\mathcal{T}_h^\Gamma$ contribute to the algebraic system.  
On $\omega_h$, we use a standard FE space of continuous functions that are piecewise polynomials of degree $m$. This bulk FE space is defined as  
\begin{equation*} 
	V_{h}^m := \{\, v \in C(\omega_h) \,:\, v|_T \in \mathcal{P}_m(T) \text{ for all } T \in \mathcal{T}_h^\Gamma \,\}.
\end{equation*}  

To understand how the TraceFEM works, let us consider the Laplace--Beltrami equations  on $\Gamma$, a model surface elliptic PDE:
\begin{equation}
	-\Delta_{\Gamma} u+u=f\quad\text{on}~~\Gamma,
	\label{LB}
\end{equation} 
and its TraceFEM formulation:
Find $u_h \in V_{h}^m$ such that
\begin{equation} \label{discr1}
	\int_{\Gh}( \unablah u_h \cdot \unablah v_h  + u_h v_h) \, \rd s + \rho_s s_h(u_h,v_h) = \int_{\Gh} f_h v_h \, \rd  s,
\end{equation}
for all $v_h \in  V_{h}^m$ with $\rho_s>0$.
Here $f_h$ denotes  an approximation of the data $f$  on $\Gamma_h$. The first and second terms in \eqref{discr1} are a FE analogue of the weak formulation of the Laplace--Beltrami problem and use only traces of $u_h,v_h$ on $\Gamma_h$, while the purpose of the third term is provide coercivity on the bulk space  $V_{h}^m$. One common choice of $ s_h(u_h,v_h)$ is so-called volume normal stabilization bilinear form~\cite{grande2018analysis} given by 
\begin{equation}\label{eq:def-normal-gradient-volume-stab}
	s_h(u_h,v_h):=  \int_{\omega_h} \bn_h \cdot \nabla u_h  \, \bn_h \cdot \nabla v_h \, \rd{x},
\end{equation}
 $\bn_h$ the normal to $\Gamma_h$ extended to $\omega_h$. For the scaling parameter $\rho_s\gtrsim h$ 
the following coercivity estimate holds:
\begin{equation} \label{coerc1}
h^{-1} \|u_h\|_{L^2(\omega_h)}^2 \lesssim \|u_h\|_{L^2(\Gamma_h)}^2 +\rho_s s_h(v_h,v_h)  \quad \text{for all}~~v_h \in  V_{h}^m.
\end{equation}
This ensures well-posedness of \eqref{discr1}. 
If additionally $\rho_s \lesssim h^{-1}$, then optimal order error analysis in $H^1$ and $L^2$ surface norms follows from approximation properties of the traces of the polynomial functions under natural assumptions on the surface and normal vector field approximations~\cite{olshanskii2009finite,reusken2015analysis}. 

The TraceFEM has been extended in several directions and applied to 
various problems, including fluid problems and systems of PDEs posed on evolving surfaces (see, e.g.,~\cite{lehrenfeld2018stabilized,olshanskii2022tangential}). Below we discuss how the method applies to discretize the model of a phase-separated vesicle given by the system~\eqref{grache-1m}--\eqref{gracke-4}.  

For this purpose, we consider generalized Taylor--Hood  velocity and pressure FE spaces on $\OGamma$:
\begin{equation}\label{TH}
 \bV_h = (V_h^{m+1})^3, \quad Q_h = V_h^m \cap L^2_0(\Gamma). 
 \end{equation}
Although  finite elements of a different order $k$ can be used for the phase-field and chemical potential, we let  $k=m$. 
The velocity space $\bV_h$ allows finite element vector fields \revt{that are} not necessarily tangential to $\Gamma$. However, a proper treatment of the tangentiality condition $\bu \cdot \bn = 0$ is critical. Enforcing it pointwise for polynomial vector functions may lead to locking (i.e., only $\bu_h = 0$ satisfies it). Following~\cite{olshanskii2018finite}, among many other works, we impose the tangential constraint weakly by adding a penalty term to the finite element formulation. Alternatively, one could enforce it using a Lagrange multiplier~\cite{gross2018trace} or employ the surface Piola transform to define $H(\text{div})$-conforming surface finite elements~\cite{brandner2022finite,lederer2020divergence}.

For a more compact presentation, we introduce the following  bilinear forms
 related to the CH  problem:
\begin{align}
& a_\mu(\mu, v) \coloneqq \int_{\Gamma} M \nabla_{\Gamma} \mu \cdot \nabla_{\Gamma}v\,ds +  \tau_\mu s_h(\mu, v), \label{eq:a_mu} \\
& a_c(c, g) \coloneqq  \epsilon\int_{\Gamma} \gradG c \cdot \gradG g \, ds + \tau_c s_h(c, g). \label{eq:a_c}
\end{align}
Forms \eqref{eq:a_mu}--\eqref{eq:a_c} are well defined for $\mu,v, c, g \in H^1(\OGamma)$.
The analysis dictates the different scaling of stabilization parameters with respect to $h$,
\[
\tau_\mu=h,\quad \tau_c=\epsilon\,h^{-1}.
\]
 In particular, we need to control the $L^\infty(\Gamma)$ norm of the FE approximation of the phase-field parameter $c$, but not the approximation of the chemical potential.

For the numerical stability, it is crucial that the computed density and viscosity stay positive, which is not automatically enforced by the numerical model.  Assuming $\rho_1\ge\rho_2$ and $\eta_1\ge\eta_2$, we consider the following cut-off functions:
\begin{align*}
\rho(c)=\left\{
\begin{array}{cc}
  \rho_2 & c\le 0, \\
  c\rho_1+(1-c)\rho_2 & c>0,
\end{array}
\right. \\
\eta(c)=\left\{
\begin{array}{cc}
  \eta_2 & c\le 0, \\
  c\eta_1+(1-c)\eta_2 & c>0.
\end{array}
\right.
\end{align*}
We further approximate $\rho(c)$  by a smooth monotone convex and uniformly positive function by letting
$\theta^2 = \frac{\rho_1-\rho_2}2\left(\tanh(c/\alpha)+1\right)$, with $\alpha=0.1$, and $\rho(c)=\int_0^c\theta^2(t)dt+\rho_2$.
The convexity  of $\rho(c)$ plays a role in the analysis.

Consider the decomposition of a vector field on $\Gamma$ into its tangential and normal components \eqref{u_T_N}
and let $\widehat{\rho}=\rho-\frac{d\rho}{d\,c}c$.
The following forms are needed or for the surface fluid equations:
\begin{align}
& a(\eta; \bu,\bv) \coloneqq \int_\Gamma 2\eta E_s( \overline{\bu}):  E_s( \overline{\bv})\, ds+\tau \int_{\Gamma}(\bn\cdot\bu) (\bn\cdot\bv) \, ds \cl
&\quad + \beta_u s_h( \bu, \bv)   \label{eq:a} \\
& c(\rho; \bw, \bu,\bv)\coloneqq \int_\Gamma\rho\bv^T(\nabla_\Gamma\overline{\bu})\bw\, ds
 +\frac12 \int_\Gamma\widehat{\rho}(\divG \overline{\bw})\overline{\bu}\cdot\overline{\bv}\, ds, \label{eq:c} \\
&b(\bu,q) =  \int_\Gamma \bu\cdot\nabla_\Gamma q \, ds, \label{eq:b} \\
&s(p,q) \coloneqq \beta_p  s_h (p,q) \label{eq:s}
\end{align}
where $\tau>0$ is a penalty parameter to enforce the tangential constraint,
$\beta_u\ge0$  and $\beta_p\ge0$ in are stabilization parameters
defined by
\begin{equation} \label{param}
\tau=h^{-2},\quad \beta_p=h, \quad \beta_u=h^{-1}.
\end{equation}

The semi-discrete finite element method for \eqref{grache-1m}--\eqref{gracke-4} then reads: Find $c_h(t)$, $\mu_h(t):(0,T] \to V^m_h$, $\bu_h(t):(0,T] \to \bV_h$, $p_h(t) :(0,T] \to  Q_h$ satisfying initial conditions and solving
\begin{align*}
	 &(\rho\partial_t\overline{\bu}_h,\bv_h)+ c(\rho; \bu_h, \bu_h,\bv_h) + a(\eta;\bu_h,\bv_h)  + b(\bv_h,p_h)  = \cl
     &\quad -(\sigma_\gamma c_h \nablaG \mu_h, \bv_h)  + M \left(\nabla_\Gamma(\theta\overline{\bu}_h)\gradG \mu_h, \theta\bv_h\right) + (\bbf_h, \bv_h), \\
      &b(\bu_h,q_h)-s(p_h,q_h)  = 0, \\
	&\left(\partial_t c_h, v_h\right) - \left(\bu_h  c_h, \nablaG v_h\right) + a_\mu(\mu_h, v_h) = 0,  \\
		&\left(  \mu_h  - \frac{1}{\epsilon} f'_0(c_h),\,g_h\right)
		- a_c(c_h, g_h) = 0,  
\end{align*}
for all $t\in (0,T]$,  $(v_h, g_h) \in V^m_h \times V^m_h$, $(\bv_h,q_h) \in \bV_h \times Q_h$. Next, we are interested in time discretization.

At time instance $t^n=n\Delta t$, with $\Delta t=\frac{T}{N}$,  $\zeta^n$ denotes the approximation of generic
variable $\zeta(t^n, \bx)$. Further, we introduce the following notation  for a first order approximation of the
 time derivative:
 \begin{equation*}
\left[\zeta\right]_t^{n} =\frac{\zeta^{n}- \zeta^{n-1}}{\Delta t}.
\end{equation*}

We introduce now a fully discrete formulation, which decouples the phase field and fluid parts of the system on each time step.
Moreover, each individual sub-problem is linear. This allows us to achieve low computational costs. 
At time step $t^{n+1}$, perform:\\[0.9ex]
\begin{itemize}
    \item[\underline{Step 1}:] Given $\bu^n_h\in \bV_h$ and $c^n_h\in V_h^m$, find $(c^{n+1}_h, \mu^{n+1}_h) \in V^m_h \times V^m_h$
such that:
\begin{align}
&\left(\left[c_h\right]_t^{n+1}, v_h\right) - \left(\bu^{n}_h  c^{n+1}_h, \nablaG v_h\right) + a_\mu(\mu_h^{n+1}, v_h) = 0,  \label{eq:CH_FE1} \\
&\left(  \mu_h^{n+1} - \frac{\gamma_c\Delta t}{\epsilon}\left[c_h\right]_{t}^{n+1}  - \frac{1}{\epsilon} f'_0(c_h^{n}),\,g_h\right)
- a_c(c_h^{n+1}, g_h) = 0,  \label{eq:CH_FE2}
  \end{align}
for all  $(v_h, g_h) \in V^m_h \times V^m_h$.
This semi-implicit time splitting is taken from \cite{shen2015decoupled}, where it is applied to the model from~\cite{Abels2012}. An alternative to this is the scalar auxiliary variable method \cite{shen2018scalar,palzhanovSAV}, which enables the construction of efficient and accurate time discretization schemes.
\item[\underline{Step 2}:] Set 
$\theta^{n+1} =\sqrt{\frac{d\rho}{d c}(c^{n+1}_h)}$.
Find $(\bu_h^{n+1}, p_h^{n+1}) \in \bV_h \times Q_h$ such that
 \begin{align}
&  (\rho^n\left[\overline{\bu}_h\right]_t^{n+1},\bv_h)+ c(\rho^{n+1}; \bu^{n}_h, \bu^{n+1}_h,\bv_h) + a(\eta^{n+1};\bu_h^{n+1},\bv_h) \cl
& \quad \quad + b(\bv_h,p_h^{n+1}) = (\bbf_h^{n+1}, \bv_h) -(\sigma_\gamma c^{n+1}_h \nablaG \mu^{n+1}_h, \bv_h) \cl
& \quad \quad + M \left((\nabla_\Gamma(\theta^{n+1}\overline{\bu}_h^{n+1}))\gradG \mu^{n+1}_h, \theta^{n+1}\bv_h\right)  \label{NSEh1} \\
 & b(\bu_h^{n+1},q_h)-s(p_h^{n+1},q_h)  = 0,  \label{NSEh2}
 \end{align}
 for all $(\bv_h^{n+1}, q_h^{n+1}) \in \bV_h \times Q_h$
\end{itemize}

The theorem below shows that  the scheme is provably stable under relatively mild restrictions. 
\begin{theorem}\label{Th:main}
\label{Th1} Assume  $h$ and $\Delta t$ satisfy $\Delta t \le c |\ln h|^{-1} \epsilon$ and  
$$h \le c |\ln h|^{-1} \min\{\Delta t, |\ln h|^{-\frac12} \epsilon |\Delta t|^{\frac12}\}$$  
for some sufficiently small $c>0$, independent of $h$, $\Delta t$, $\epsilon$ and position of $\Gamma$ in the background mesh. Then, the solution to \eqref{eq:CH_FE1}--\eqref{NSEh2} satisfies
 \begin{multline}
 \int_\Gamma\left(\rho^{N}|\overline{\bu}_h^{N}|^2 +\frac{\sigma_\gamma}{\epsilon}f_0(c^{N}_h)\right) ds  + a_c(c^{N}_h,c^{N}_h) \\ + \sum_{n=1}^{N}\Delta t \left( a(\eta^{n}; \bu^{n}_h,\bu^{n}_h) + a_\mu(\mu^{n}, \mu^{n})+ s_h(p_h^{n},p_h^{n})\right)
 \le
 K,  \label{stab_est}
\end{multline}
for all $N=1,2,\dots$, with $K=\int_\Gamma\left(\rho^{0}|\bu^{0}_h|^2 +\frac{\sigma_\gamma}{\epsilon}f_0(c^{0}_h)\right) ds + a_c(c^{0}_h,c^{0}_h)$.
\end{theorem}
For the proof, we refer to \cite{palzhanov2021decoupled}. However, we note
that the results of extensive numerical experiments in~\cite{palzhanov2021decoupled} do not show that any restriction on the discretization parameters is required in practice.

\section{Validation against experimental data} \label{sec:Exp}

\subsection{NSCH model}
\label{sec:NSCH_val}  

We consider two membrane compositions exhibiting distinct and nearly opposite  phase behaviors: one with a majority $L_o$ phase and another with a minority $L_o$ phase.
These two membrane compositions are DOPC:DPPC:Chol at molar ratio of 1:1:15\%, in which the $L_o$ phase is expected 
to occupy about 29\% of the membrane surface at 25\textdegree{}C, and 1:2:25\%, in which the $L_o$ phase would occupy about 
71\% of the membrane area at 15\textdegree{}C. 
These indicated area fractions, which will be denoted with
$a_D$, are calculated using an approach 
that relies on the composition of each phase, determined with the phase diagram tie-lines (see example in Fig.~\ref{fig:PD}) and the molecular area of the lipid components. See \cite{zhiliakov2021experimental} for more details.
The CH model from Sec.~\ref{sec:CH}, as well as continuum based models applied in other  studies~\cite{Wang2008,lowengrub2009phase,sohn2010dynamics,Li_et_al2012,Funkhouser_et_al2014}, would predict nearly the same evolution
of the domain ripening process for these two compositions since it does not account for in-membrane viscous and transport effects. However, 
we shall see that the experimental data 
reveal different evolutions that can be better  captured by the more complex 
NSCH model described in Sec.~\ref{sec:NSCH}.

\revo{ 
A modified version of electroformation \cite{kang2013simple,zhiliakov2021experimental} was employed
to form GUVs. 
A mixture of DOPC, DPPC, and Chol, plus 0.3 mol\% 
Rho-PE and 0.5\% NAP to enable fluorescence microscopy,
was prepared in chloroform and used to produce a thin lipid film in a flask using a rotary evaporator. The lipid film was then rehydrated and the resultant milky suspension was tip-sonicated to produce a clear suspension of small vesicles. GUVs were harvested from the electroformation chamber and placed on a clean microscope glass slide for imaging. Prior to imaging, 
the sample was heated on a hot plate to $\sim$60\textdegree{}C for 5 min and then placed on the microscope stage 
where it gradually cooled down to the room temperature. The image collection time was recorded with time zero 
considered as when the sample was removed from the hot plate. Epi-fluorescence microscopy was used for the 
initial assessment of GUVs and their lipid domains while confocal microscopy was used to further assess the domains 
on GUVs and quantify their size.
}
The reader interested in \revo{more} details \revo{about} the experimental set-up 
is referred to \cite{wang2022lipid}. 

For the numerical results, we considered a mesh with 225822 active degrees of freedom (193086 for $\bu_h$
and 10912 for $p_h$, $c_h$, and $\mu_h$).
The time step $\Delta t$ adaptively varies
from  $\Delta t=$4$\times 10^{-6}$ s during the fast initial phase of spinodal decomposition to about $\Delta t=$8$\times 10^{-4}$ s
during the later slow phase of lipid domain coarsening and growth, and up to $\Delta t= 4$ s when the process is close to equilibrium (around 4000 s). In order to model an initially homogenous liposome, the surface fraction $c_0$
is defined  as a realization of Bernoulli random variable~$c_\text{rand} \sim \text{Bernoulli}({a_D})$
with mean value ${a_D}$, 
i.e.,~we set:
\begin{equation}\label{raftIC}
	c_0 \coloneqq c_\text{rand}(\bx)\quad\text{for active mesh nodes $\bx$}.
\end{equation}
As mentioned above,
we set ${a_D}=0.71$ for the 1:2:25\% composition and ${a_D}=0.29$ for the 1:1:15\% composition. For each composition, we ran 10 numerical simulations.

\begin{figure*}[htb!]
\begin{center}
\vskip .5cm
\begin{overpic}[width=.18\textwidth,grid=false]{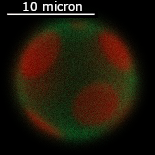}
\put(33,102){\small{$t = 102$}}
\end{overpic}~~
\begin{overpic}[width=.18\textwidth,grid=false]{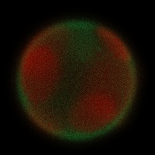}
\put(33,102){\small{$t = 145$}}
\put(48,117){\small{Composition 1:2:25\%}}
\end{overpic}~~
\begin{overpic}[width=.18\textwidth,grid=false]{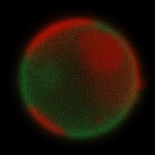}
\put(33,102){\small{$t = 408$}} 
\end{overpic}~~
\begin{overpic}[width=.18\textwidth,grid=false]{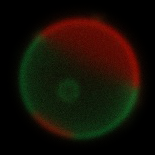}
\put(28,102){\small{$t = 1030$}}
\end{overpic}
\\
\begin{overpic}[width=.2\textwidth,grid=false]{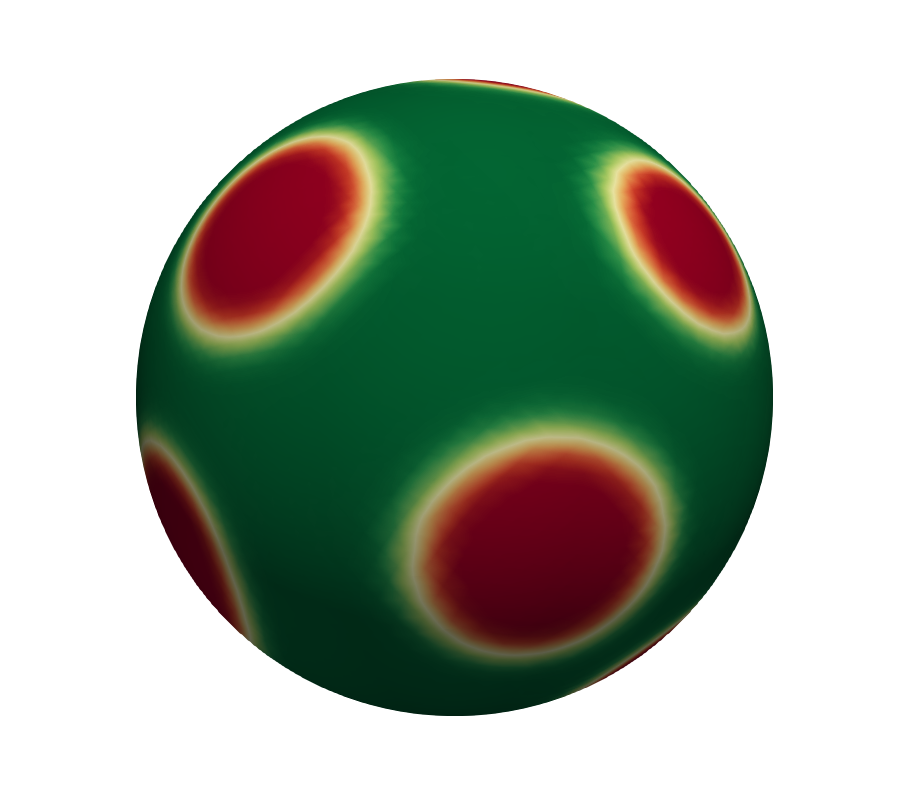}
\end{overpic}
\begin{overpic}[width=.2\textwidth,grid=false]{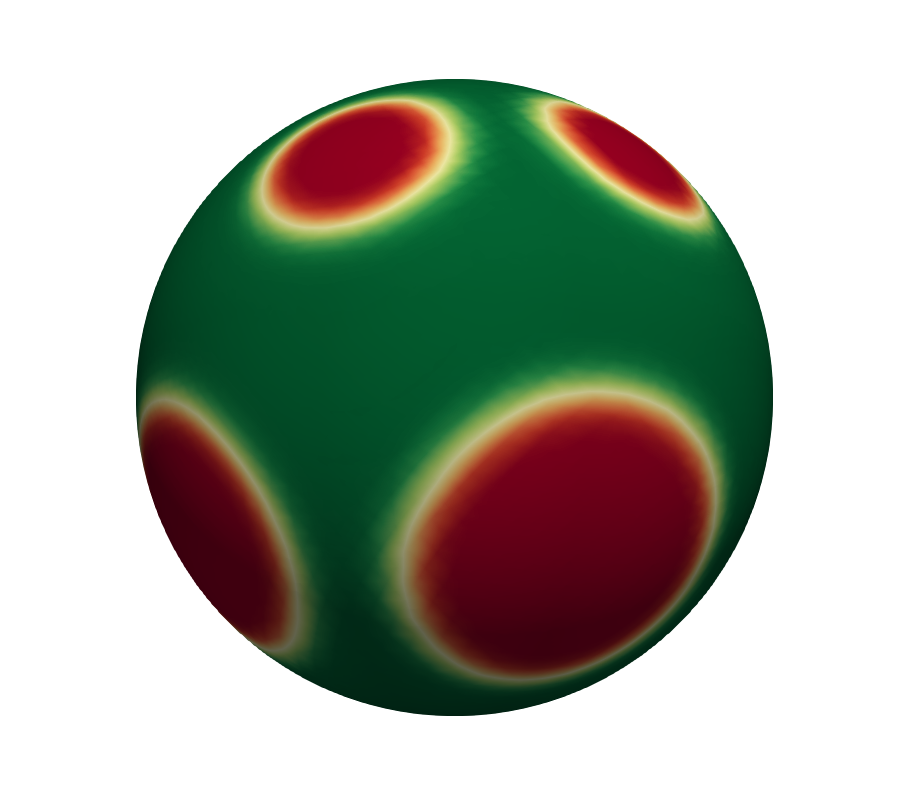}
\end{overpic}
\begin{overpic}[width=.2\textwidth,grid=false]{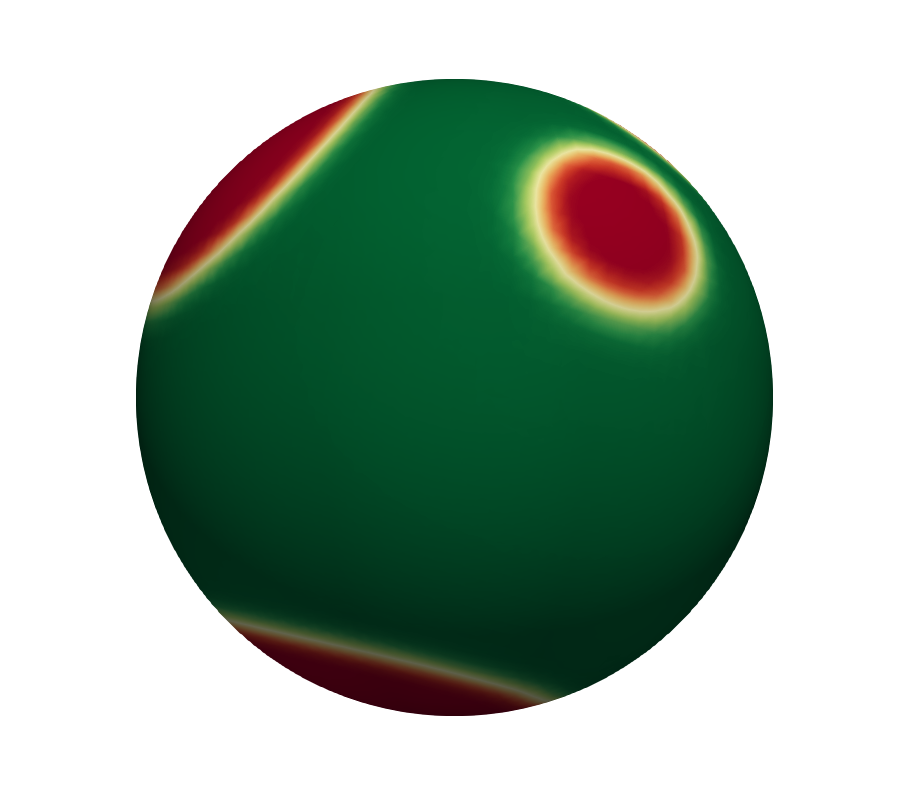}
\end{overpic}
\begin{overpic}[width=.2\textwidth,grid=false]{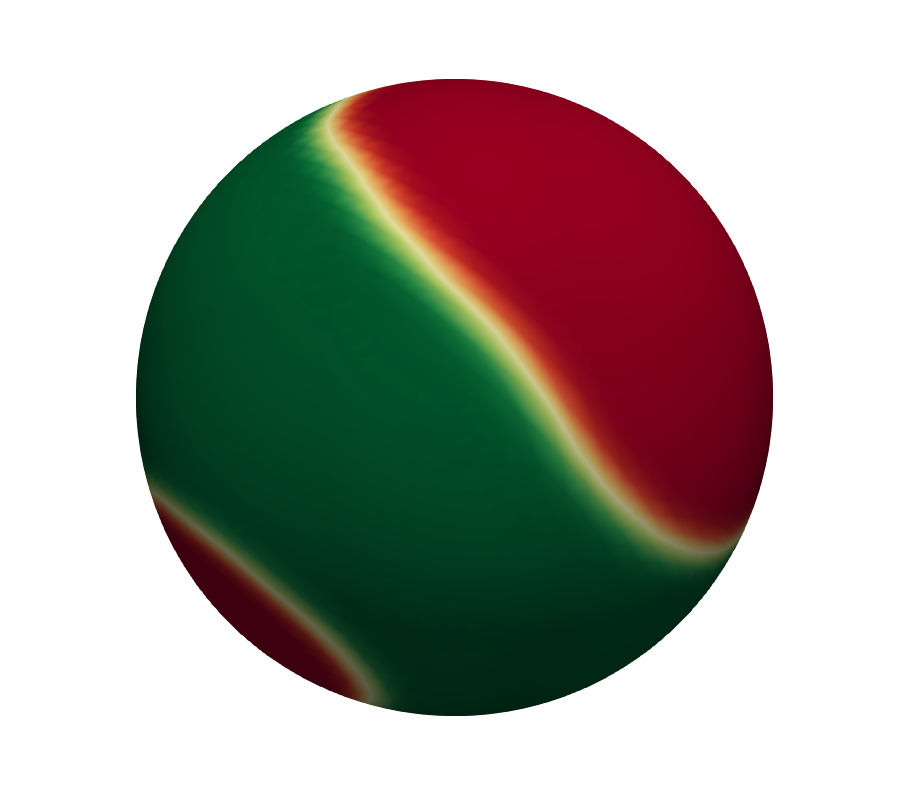}
\end{overpic}
\\
\vskip .8cm
\begin{overpic}[width=.18\textwidth,grid=false]{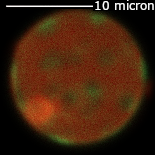}
\put(30,102){\small{$t = 73$}}
\end{overpic}~~
\begin{overpic}[width=.18\textwidth,grid=false]{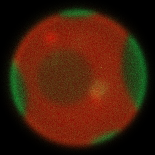}
\put(33,102){\small{$t = 166$}}
\put(48,117){\small{Composition 1:1:15\%}}
\end{overpic}~~
\begin{overpic}[width=.18\textwidth,grid=false]{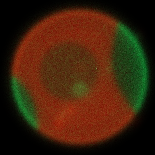}
\put(33,102){\small{$t = 244$}}
\end{overpic}~~
\begin{overpic}[width=.18\textwidth,grid=false]{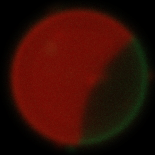}
\put(33,102){\small{$t = 322$}}
\end{overpic}
\\
\vskip .2cm
\begin{overpic}[width=.17\textwidth,grid=false]{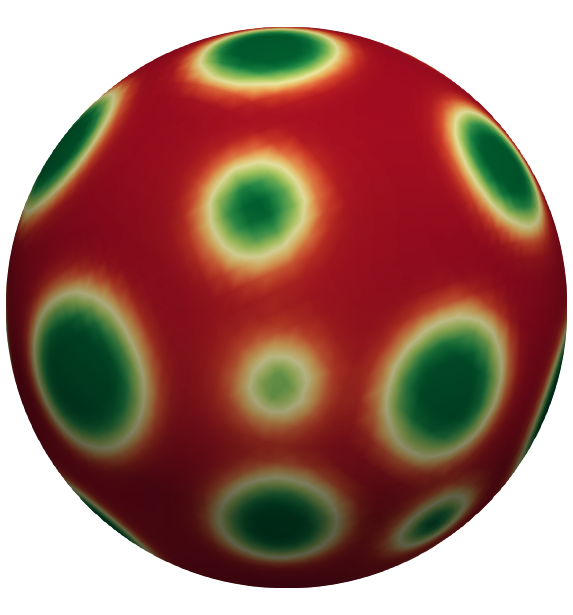}
\end{overpic}~~~
\begin{overpic}[width=.17\textwidth,grid=false]{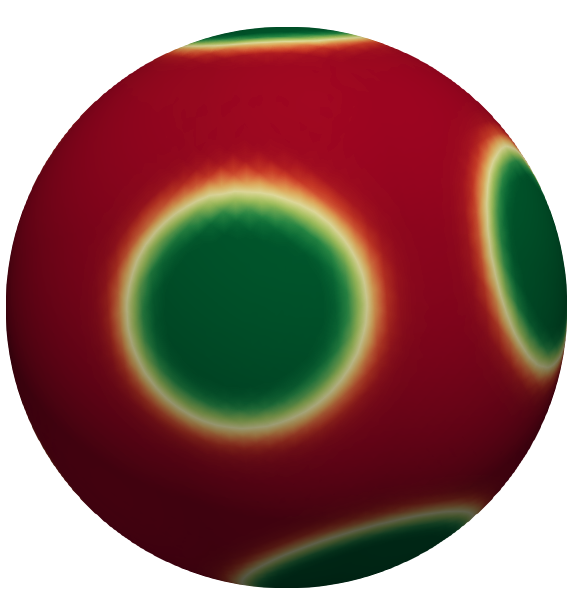}
\end{overpic}~~~
\begin{overpic}[width=.17\textwidth,grid=false]{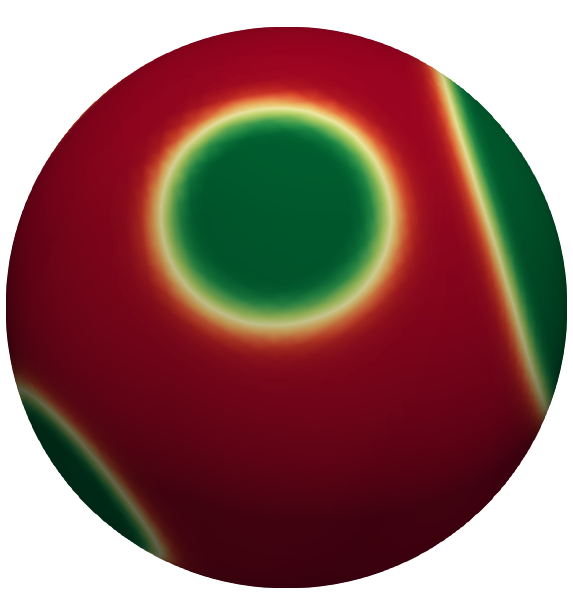}
\end{overpic}~~~
\begin{overpic}[width=.17\textwidth,grid=false]{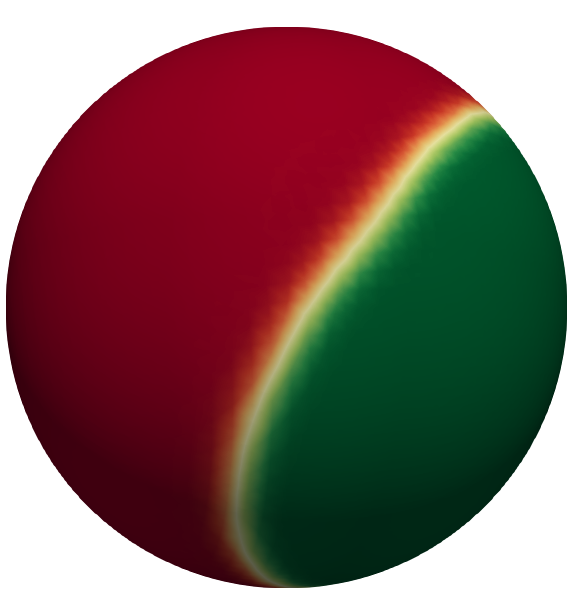}
\end{overpic}
\end{center}
\caption{Qualitative comparison for 1:2:25\% (top two rows) and 1:1:15\% (bottom two rows): epi-fluorescence microscopy images (with black background)
and numerical results (with white background) at four different times.}\label{fig:qualitative}
\end{figure*}

Let us start with a qualitative 
comparison between images acquired with
epi-fluorescence microscopy and images obtained from post-processing the numerical results.
Fig.~\ref{fig:qualitative} shows such comparison for compositions 
1:2:25\% and 1:1:15\%. 
Small quantities of dye have been added to the GUVs so that the $L_o$ 
phase is colored in green, while the red regions represent 
the $L_d$ phase. The colors for the numerical results have been chosen accordingly.
Overall, from Fig.~\ref{fig:qualitative} we see an excellent
qualitative agreement between experiments and simulations.

\begin{figure*}[htb!]
\centering
		\begin{overpic}[width=.42\textwidth, grid=false]{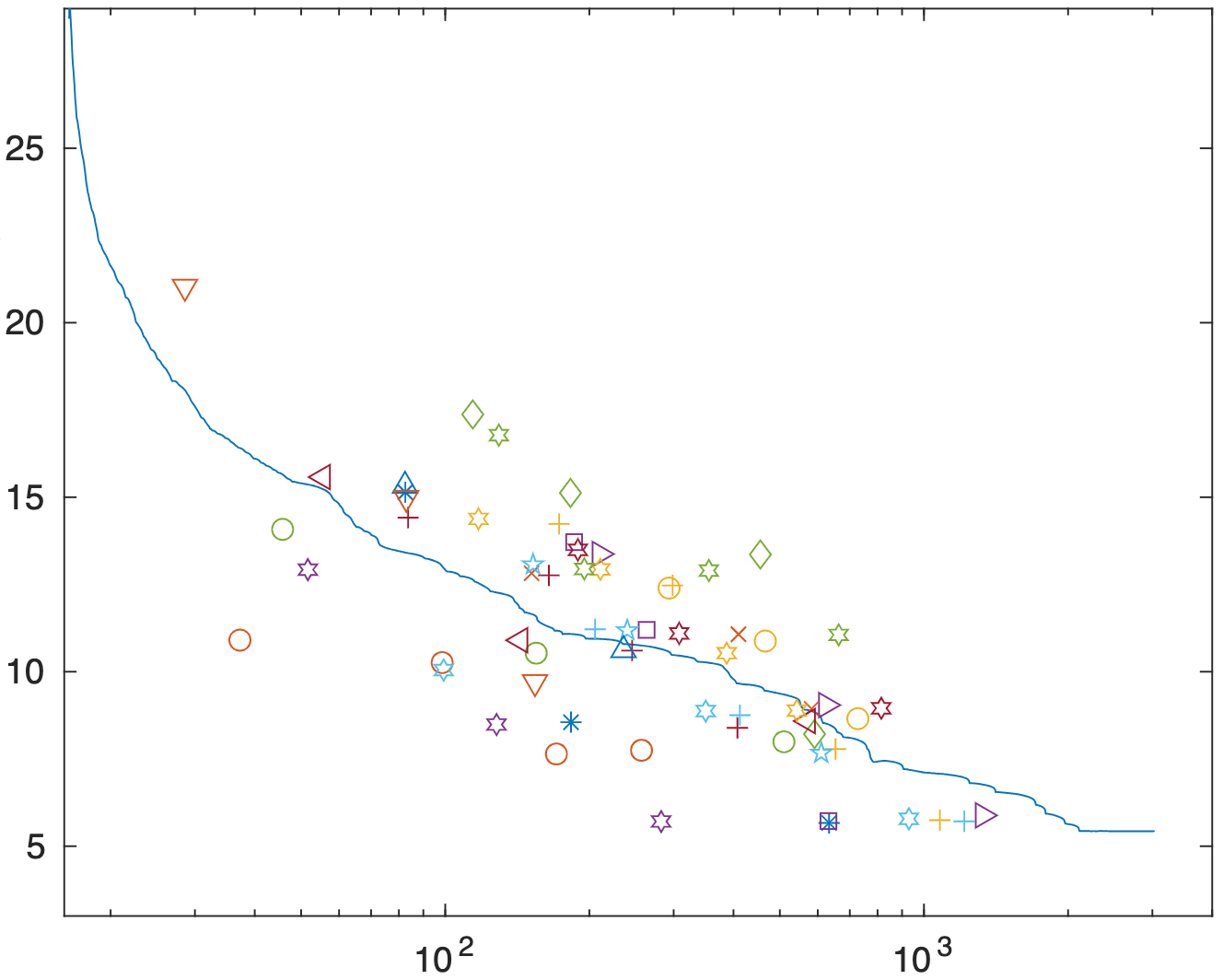} 
        \put(47,-3){\footnotesize{time \revo{(s)}}}
        \put(-8,45){\footnotesize{$p_{\text{ld}}$}}
        \end{overpic}
        \quad \quad
        \begin{overpic}[width=.42\textwidth, grid=false]{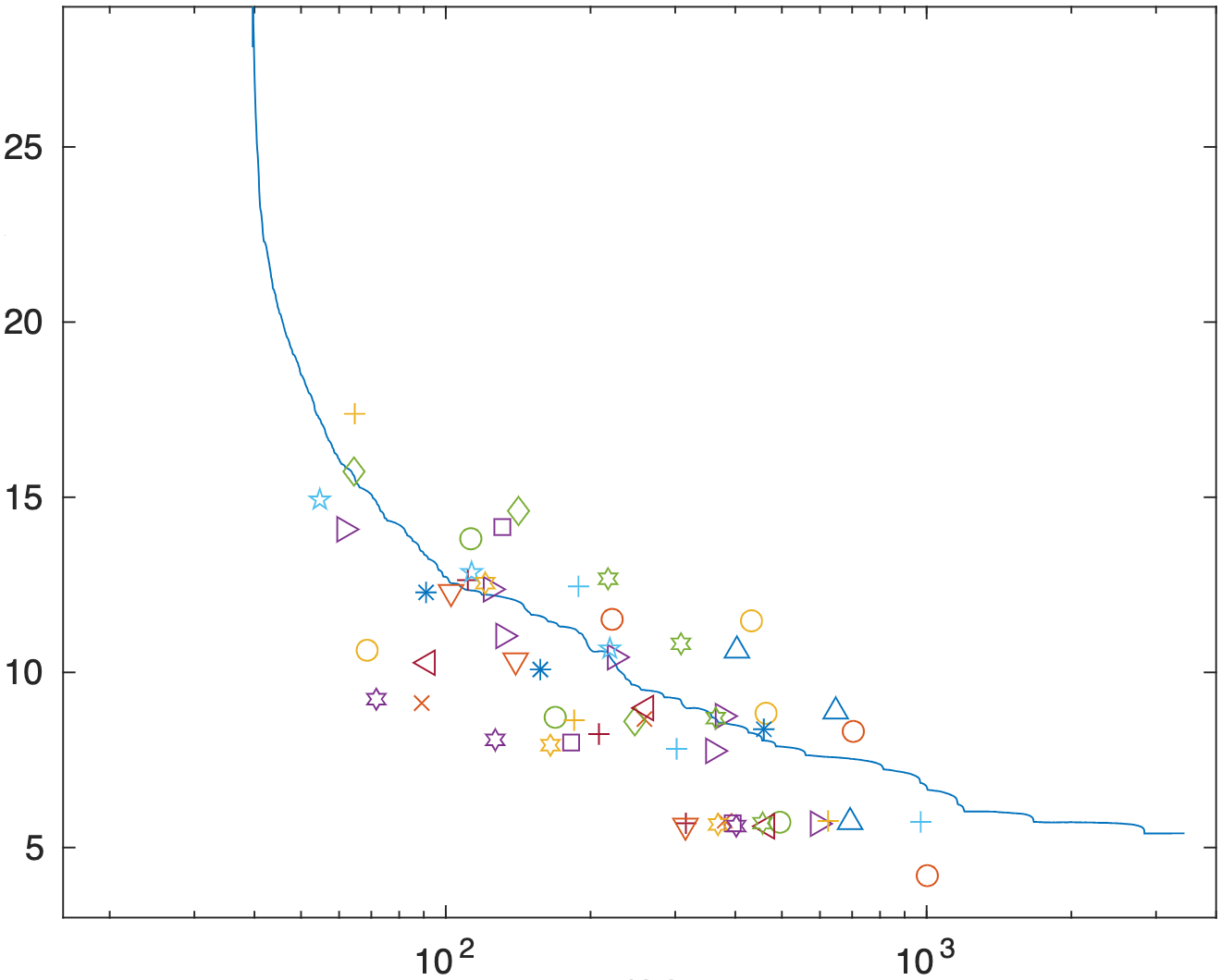}
        \put(47,-3){\footnotesize{time \revo{(s)}}}
        \put(-8,45){\footnotesize{$p_{\text{ld}}$}}
        \end{overpic}
	\caption{Total lipid domain perimeter $p_{\text{ld}}$ in \SI{}{\micro\metre} over time for composition 1:2:25\% (left)
	 and 1:1:15\% (right): numerical results average (solid line) and experimental data (markers).
     \revo{Different markers correspond to different GUVs analyzed experimentally.}}
	\label{fig:perimeter}		
\end{figure*}

\begin{figure*}[htb!]
\centering
\begin{overpic}[width=.42\textwidth, grid=false]{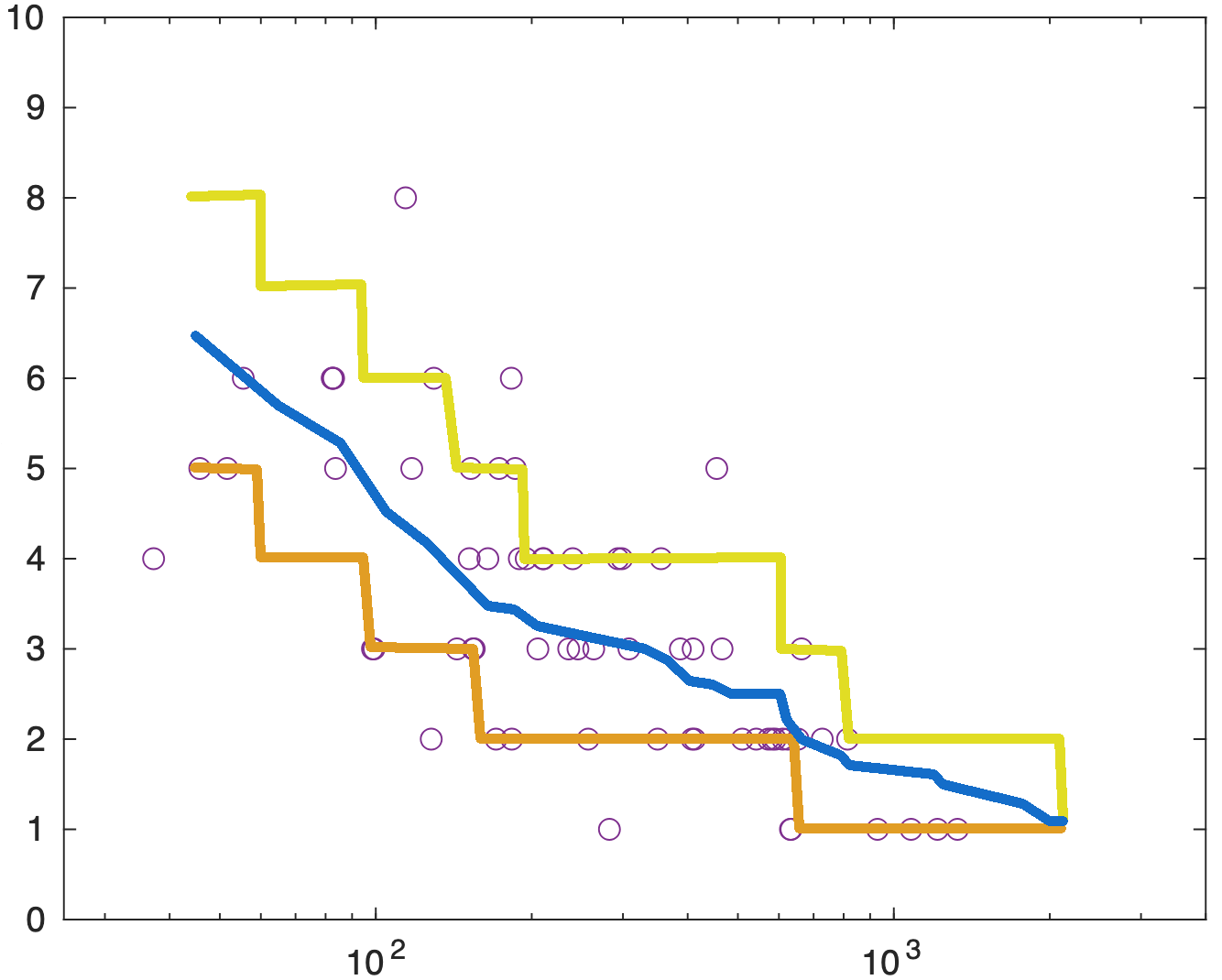} 
        \put(42,-3){\footnotesize{time (s)}}
        \put(-3,43){\makebox(0,0){\rotatebox{90}{\footnotesize{number of lipid domains}}}}
        \end{overpic}
        \quad \quad
        \begin{overpic}[width=.42\textwidth, grid=false]{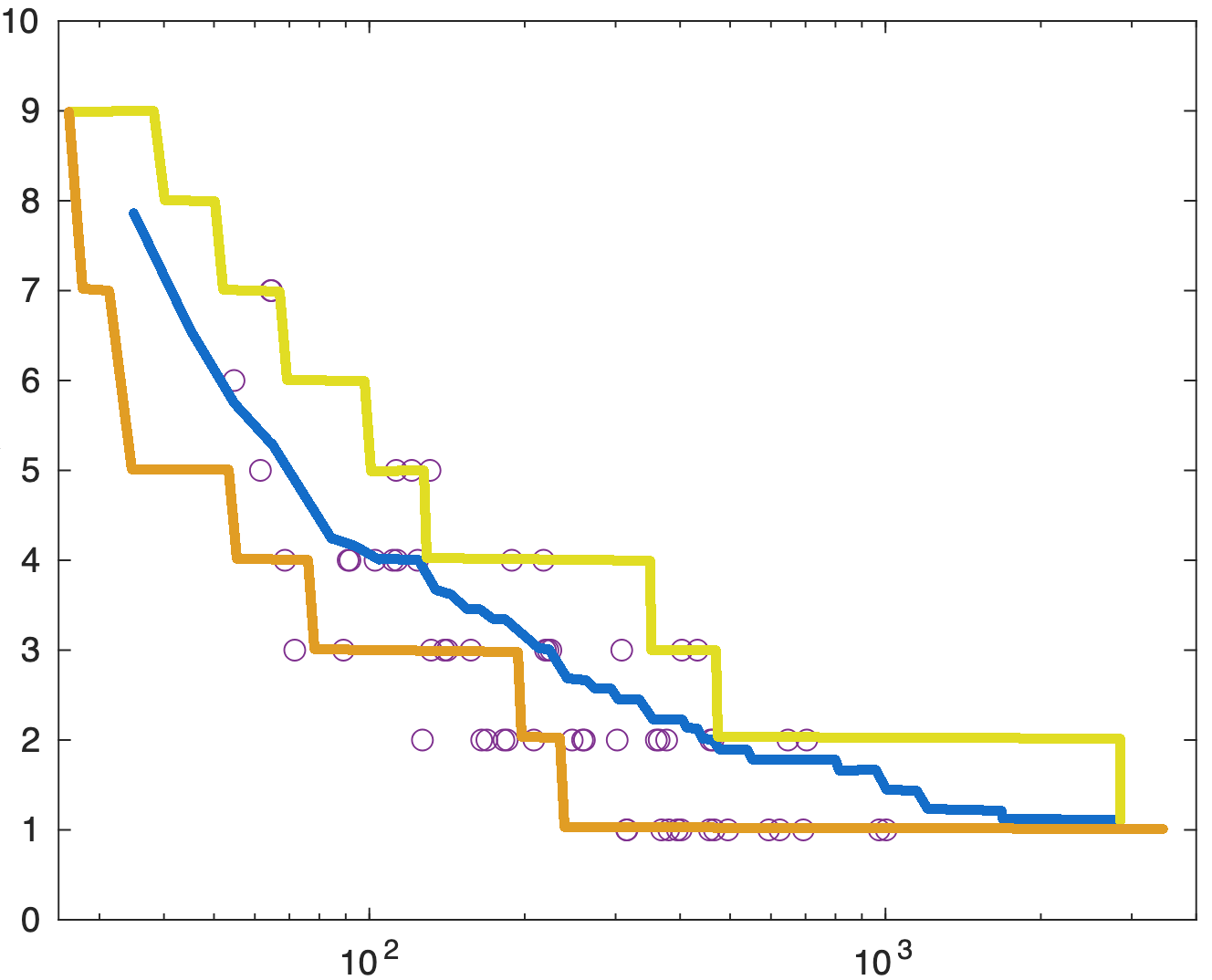}
        \put(42,-3){\footnotesize{time (s)}}
        \put(-3,43){\makebox(0,0){\rotatebox{90}{\footnotesize{number of lipid domains}}}}
        \end{overpic}
	\caption{Total number of lipid domains
	over time for composition 1:2:25\% (left) and 1:1:15\% (right): numerical results average (solid blue line),
	 minimum and maximum values found numerically (solid orange and yellow lines, respectively), and experimental data (circles).}
	\label{fig:Raft}		
\end{figure*}

For a quantitative comparison between experimental
data and numerical results, we tracked experimentally and 
numerically the total lipid domain perimeter and the total number of 
lipid domains.
We remark that numerically the total lipid domain perimeter $p_{\text{ld}}$ is computed as
\begin{equation}\label{perimeter}
p_{\text{ld}}(t_n):=2\pi\int_{\Gamma_h}\epsilon|\nabla_\Gamma c_h(\bx, t_n)|^2 ds.
\end{equation}
Fig.~\ref{fig:perimeter} reports all the experimental
measurements with markers (a different marker for each GUV) and the average of the computed total lipid domain perimeter
from all the simulations with a solid line for compositions 1:2:25\% and 1:1:15\%. In both cases,
the average of the computed total lipid domain perimeters falls within the cloud of experimental measurements.
We note that time starts at \SI{40}{\second} because no lipid domains were observed 
before \SI{40}{\second},  presumably due to the small size of domains that could not be resolved under fluorescence microscopy.
Next, we compare the total number of lipid domains on a GUV over time.
Fig.~\ref{fig:Raft} shows the experimentally measured and numerically computed data for both compositions under consideration.
The measurements are reported with a circle, while for the simulations we reported three solid lines corresponding
to the numerical results average, minimum, and maximum number of lipid domains found in the simulations.
We see that the vast majority of the experimental data 
falls within the computed extrema.

\begin{figure*}[htb!]
		\centering
        \begin{overpic}[width=.42\textwidth, grid=false]{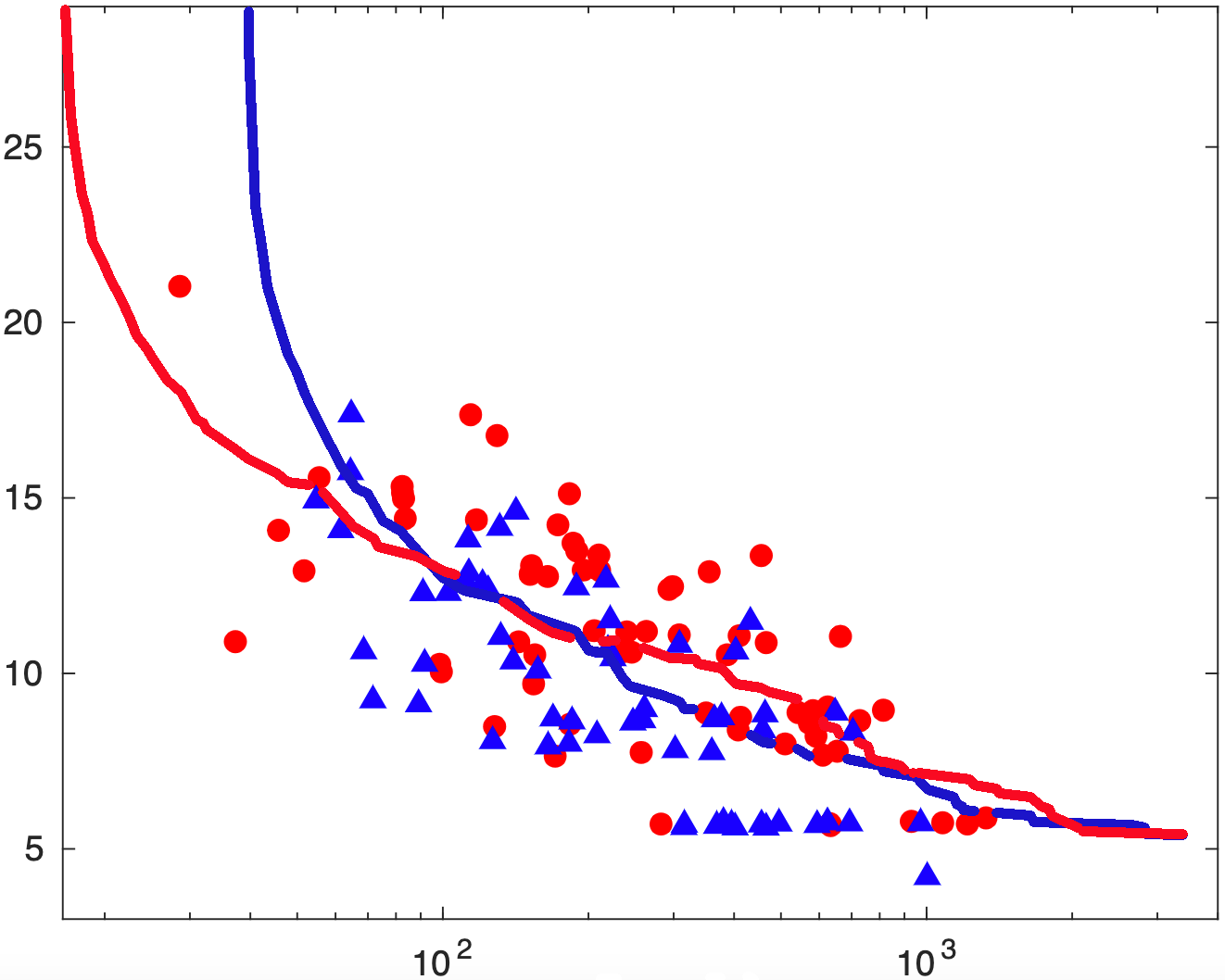} 
        \put(42,-3){\footnotesize{time (s)}}
        \put(-8,45){\footnotesize{$p_{\text{ld}}$}}
        \end{overpic}  \quad \quad
        \begin{overpic}[width=.42\textwidth, grid=false]{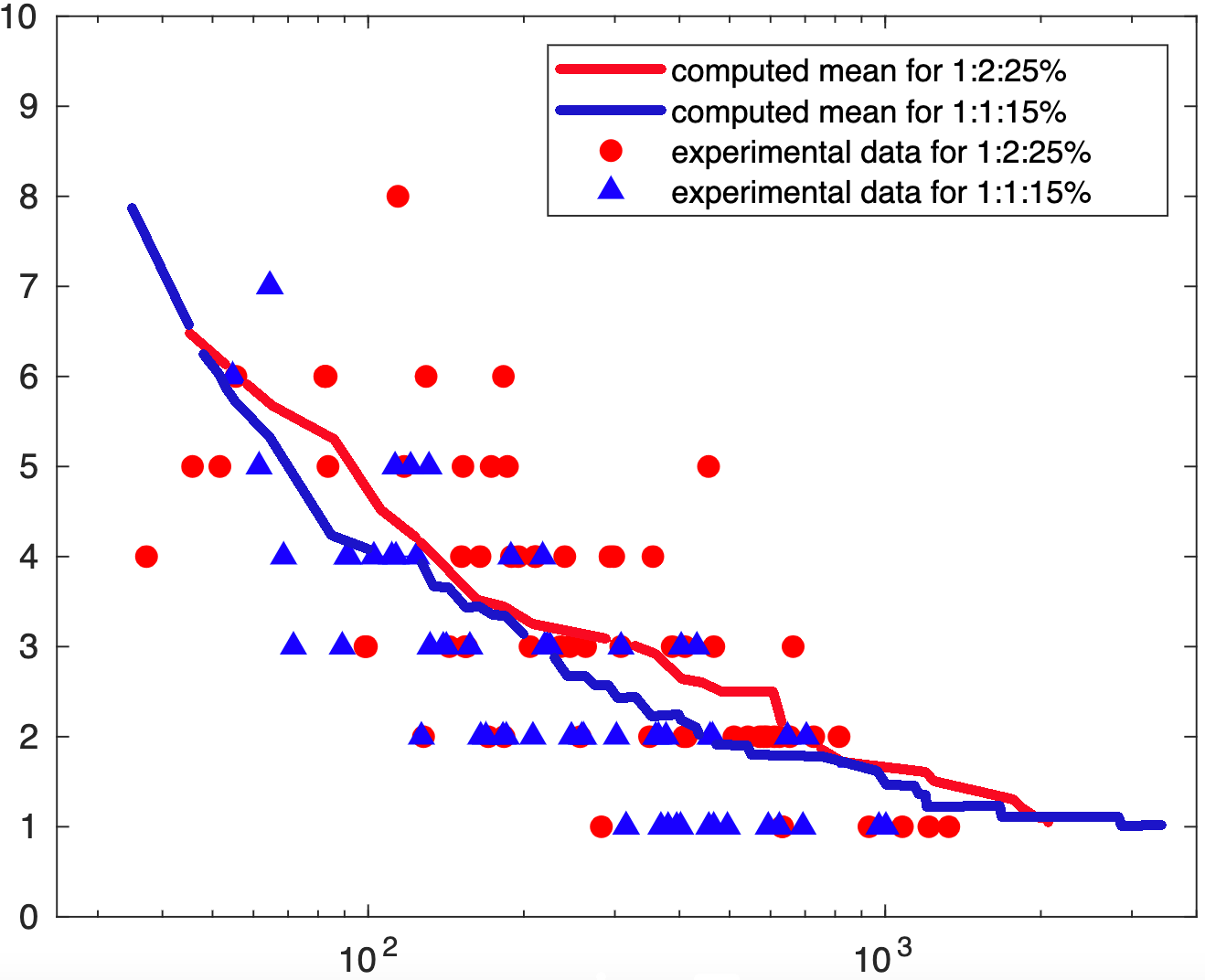} 
        \put(42,-3){\footnotesize{time (s)}}
        \put(-3,43){\makebox(0,0){\rotatebox{90}{\footnotesize{number of lipid domains}}}}
        \end{overpic}
	\caption{Experimental data for composition 1:2:25\% (red dots) and 1:1:15\% (blue triangles)
	with the corresponding computed means (solid line with corresponding color) for the total lipid domain perimeter (left) and
	total number of lipid domains (right). The legend in the subfigure on the right is common to both subfigures.
    }
	\label{fig:superimposed}		
\end{figure*}

We compare the evolution of the domain
ripening process 
for our two membrane compositions in terms of
total lipid domain perimeter and total number of lipid domains in Fig.~\ref{fig:superimposed}.
We observe in average faster dynamics towards the equilibrium state (i.e., one domain of the minority
phase within a background of the majority phase) for composition 1:1:15\%. 
This is correctly captured by the NSCH model described in Sec.~\ref{sec:NSCH}.
Indeed, we see that the solid blue curve (corresponding to the computed mean for composition 1:1:15\%) lies below the red curve
(corresponding to the computed mean for composition 1:2:25\%) for the majority of the time interval under consideration
in both graphs in Fig.~\ref{fig:superimposed}. 
The CH model from Sec.~\ref{sec:CH} would predict 
the same evolution for both composition and thus it would be unsuited
to reproduce the experimental data. 

\subsection{NSCH model with an external force}\label{sec:charge_val}

Let us turn to phase separation in  DOTAP:DOPC:DPPC:Chol mixture to modulate surface density of DOTAP (a cationic lipid) on liposomes and hence their fusogenicity. Recall that 
the idea is to concentrate DOTAP, through phase separation,  into small patches on the liposome's surface with the goal of enhancing the liposome's fusogenicity without the need for high DOTAP concentrations, which are toxic in vivo. See Sec.~\ref{s1}.
 
We consider SUVs of three different phase-separating compositions containing DOTAP (referred to as patchy liposomes - PAT). The fusogenicity of these SUVs 
into GUVs, acting as model target membranes,
is compared to the fusogenicity of homogeneous SUVs, i.e., SUVs that do not undergo phase separation. For the phase-separating SUVs, lipid composition 
DOPC:DPPC:Chol with three different molar ratios are selected (see Table \ref{tab:comp}), 
in which 15\% of DOPC is replaced with DOTAP.
Given the similarity between DOTAP and DOPC, 
DOPC mostly partitions into the $L_d$ phase \cite{wang2024}. 
With the same DOTAP content, composition PAT3 is expected to have the highest surface density of DOTAP in $L_d$ phase
because it has the largest $a_D$, 
and composition PAT1 is expected to have the lowest density of DOTAP in its $L_d$ phase because it has the smallest $a_D$. 

\begin{table}[htb]
\begin{center}
 \begin{tabular}{ | c |  c |  c |  c |  c |}
\hline
Composition  & DOPC & DPPC & Chol & $a_D$ \\
\hline
Homo & 99.4\% & 0\% & 0\% &0\%\\
\hline
PAT1 & 59.4\% & 20\% & 20\% & 10.8\% (15\textdegree{}C) \\
\hline
PAT2 & 41.9\% & 42.5\% & 15\% & 34.57\% (17.5\textdegree{}C) \\
\hline
PAT3 & 24.4\% & 50\% & 25\% & 70.37\% (15\textdegree{}C)\\
\hline
\end{tabular}
\caption{Lipid composition for the examined liposomes.
}\label{tab:comp}
\end{center}
\end{table}

\begin{figure}[htb!]
	\centering
	\includegraphics[width = .48\textwidth]{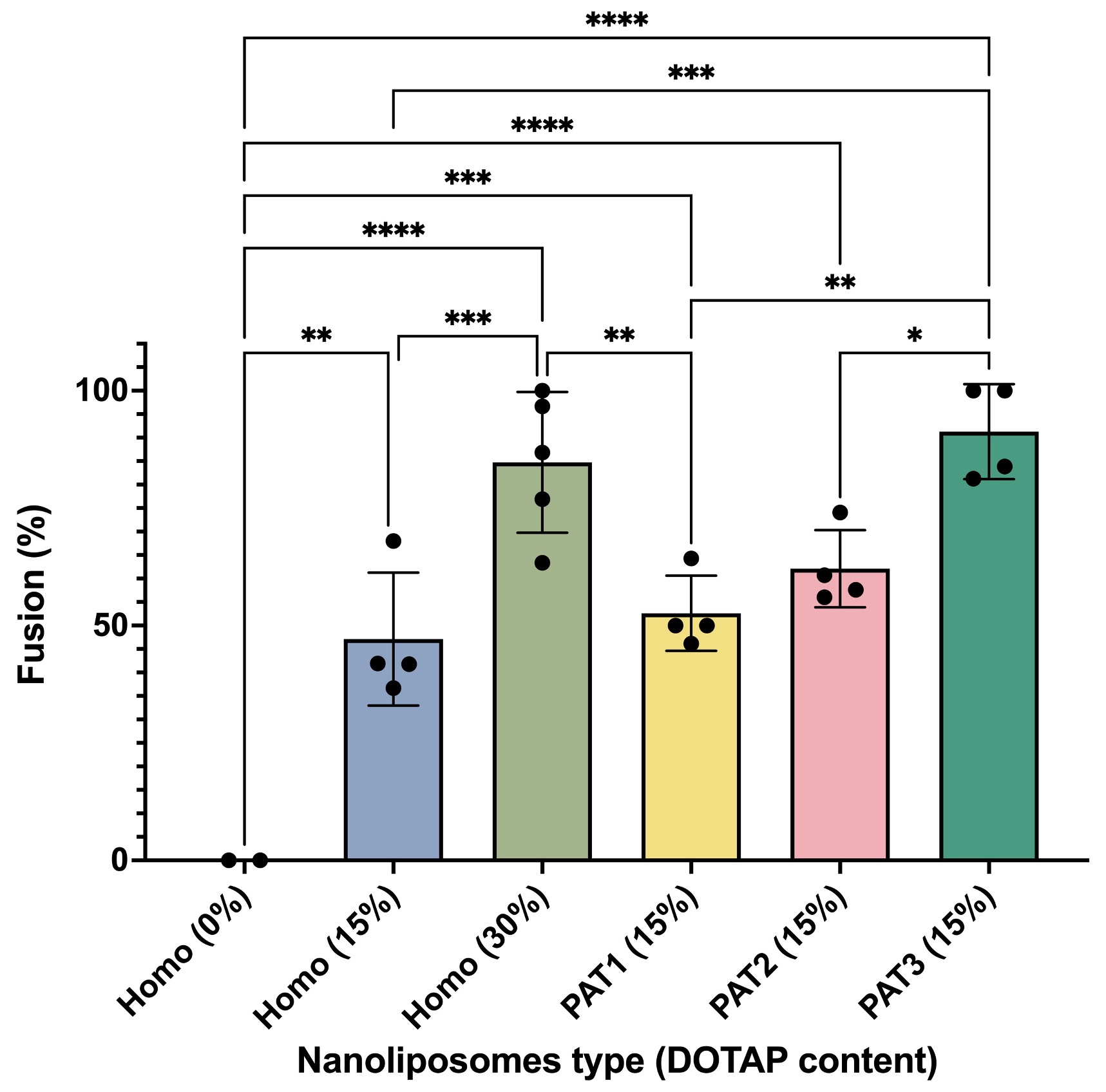}
	\caption{Fraction of GUVs that showed fusion after \SI{10}{\minute} incubation with SUVs of different lipid compositions from \cite{wang2024}.
The data were statistically analyzed using one-way ANOVA and $^*$: p value$<$0.05, $^{**}$: p value$<$0.01, $^{***}$: p values$<$0.001, $^{****}$: p value$<$0.0001.
}
	\label{fig:fusion_frac}
\end{figure}

In the experiments reported in \cite{wang2024}, homogeneous and phase-separated
SUVs were incubated with GUVs of DOPC composition at 37\textdegree{}C for 10 min. 
After the incubation, samples were imaged with confocal microscopy to evaluate the level of fusion of SUVs 
(labeled with red fluorescence) into GUVs (labeled with green fluorescence). In case of homogenous SUVs with no DOTAP, 
the GUVs exhibited only green fluorescence indicating no significant fusion. 
Upon increasing DOTAP concentration in homogenous SUVs, a mixture of both 
red and green fluorescence on GUVs was obtained, suggesting some level of fusion. 
To quantify the level of fusion in the experiments, the fraction of 
GUVs that showed fusion upon incubation with SUVs
was measured. Fig.~\ref{fig:fusion_frac} reports the results: higher DOTAP concentration result in higher level of fusion and PAT3 composition, 
with highest DOTAP density in $L_d$ phase, shows the highest level of fusion. In particular, note that 
phase-separating SUVs of PAT3 composition (with 15\% DOTAP), 
led to much stronger red fluorescence signal in GUV membranes compared to that in case of homogeneous liposomes with 
15\% DOTAP, and was comparable to that of homogeneous SUVs with 30\% DOTAP. 

To reproduce the experiments, one phase-separated SUV
in equilibrium state, i.e., one patch of the minority phase
against the background of the majority phase is exposed to one GUV. See Fig.~\ref{fig:force}. The dynamics of the
phases and the surface flows are simulated with model 
\eqref{grache-1m}-\eqref{gracke-4}, with $\bbf = \bbf_e$, as explained in Sec.~\ref{sec:ex_forces}. For additional clarity,
let us detail  the calculation of $\bbf_e$ for a PAT3 SUV, 
which has $a_D = 70.37\%$ (see Tab.~\ref{tab:comp}), i.e., about 70\% of the surface of the SUV
is covered by the $L_o$ phase (red in Fig.~\ref{fig:force}). For composition PAT3, the concentration of DOTAP 
in the $L_d$ phase (blue in Fig.~\ref{fig:force}) is 41.8\% (see Tab.~4 in \cite{wang2024}), corresponding to 67.15\% of the total DOTAP in the SUV.  
So, we uniformly distribute 67.15\% of the total charge density, and hence \revt{the} force, to the $L_d$ phase. 
For the simulations, we adopt the same mesh and time steps used for the results in Sec.~\ref{sec:NSCH_val}.

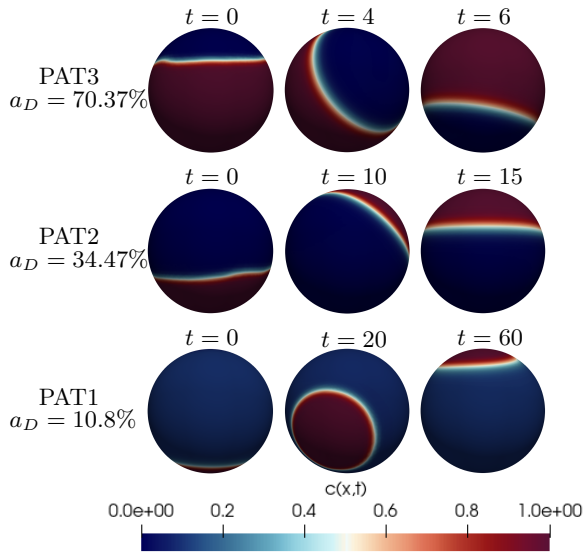
\begin{figure*}[htb!]
	\begin{center}
        \begin{overpic}[width=.13\textwidth,grid=false]{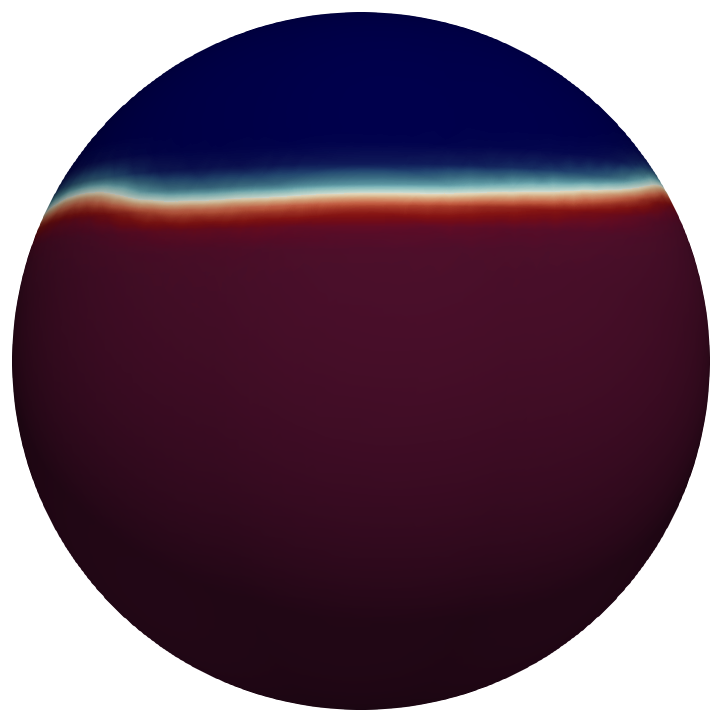}
			\put(30,102){\small{$t = 0$}}
			\put(-85,55){\small{PAT3}}
			\put(-108,37){\small{$a_D = 70.37$\%}}
		\end{overpic}
		\begin{overpic}[width=.13\textwidth,grid=false]{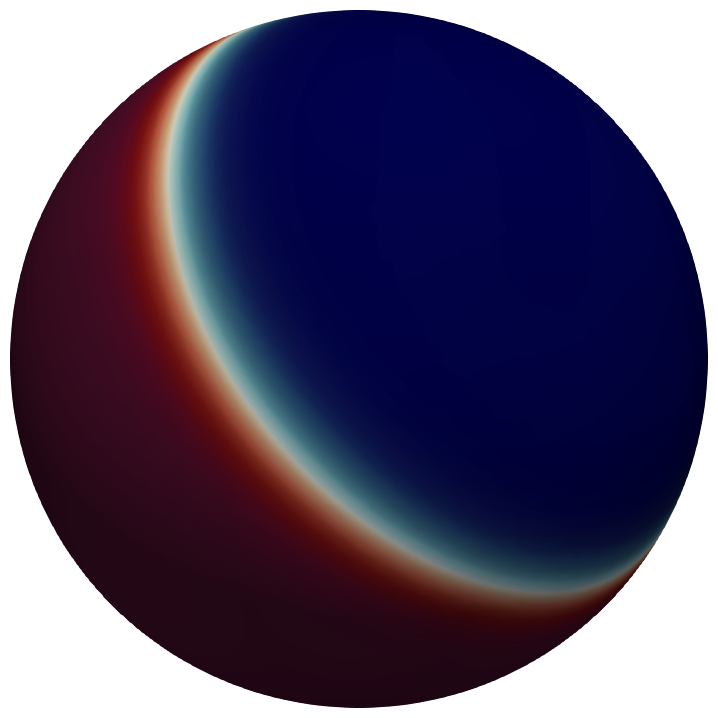}
			\put(30,102){\small{$t = 4$}}
		\end{overpic}
		\begin{overpic}[width=.13\textwidth,grid=false]{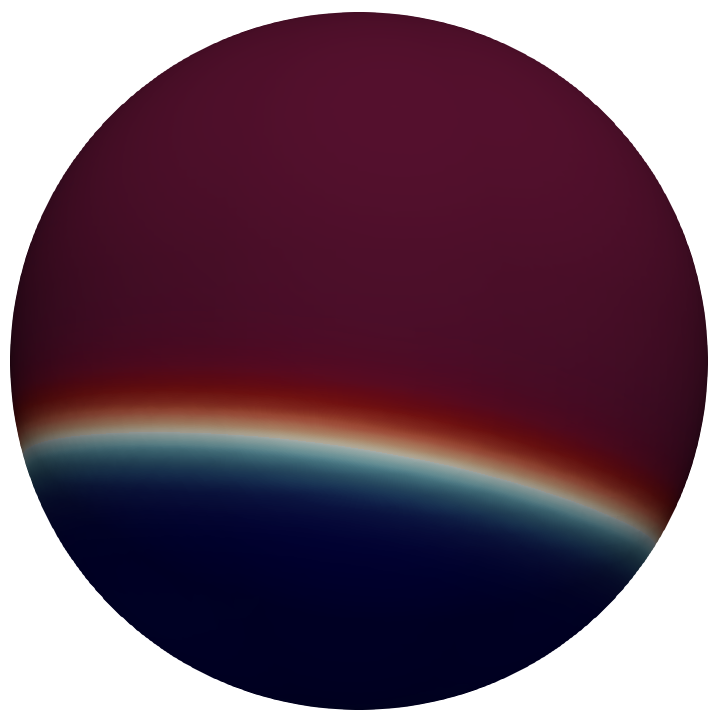}
			\put(30,102){\small{$t = 6$}}
		\end{overpic}\\
		
		\vskip 12pt
        \begin{overpic}[width=.13\textwidth,grid=false]{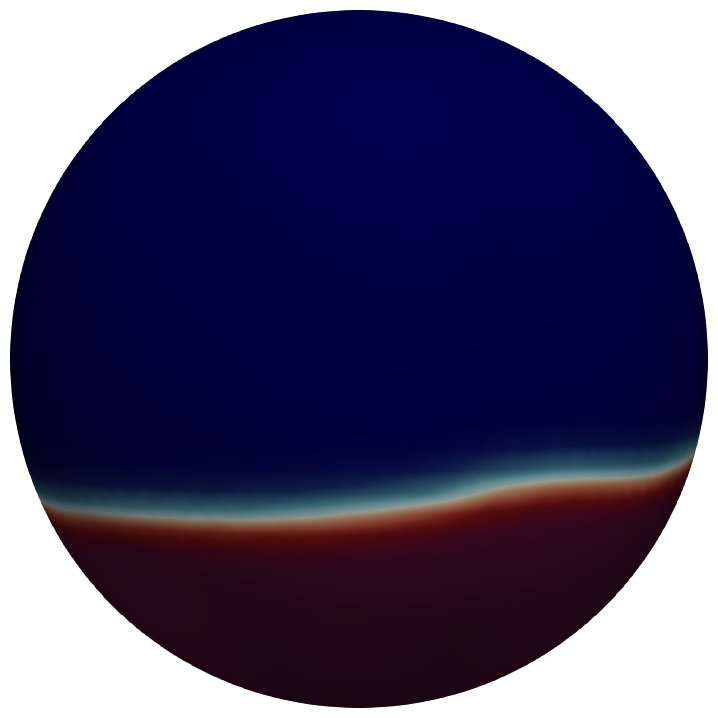}
			\put(30,102){\small{$t = 0$}}
			\put(-85,55){\small{PAT2}}
			\put(-108,37){\small{$a_D = 34.47$\%}}
		\end{overpic}
		\begin{overpic}[width=.13\textwidth,grid=false]{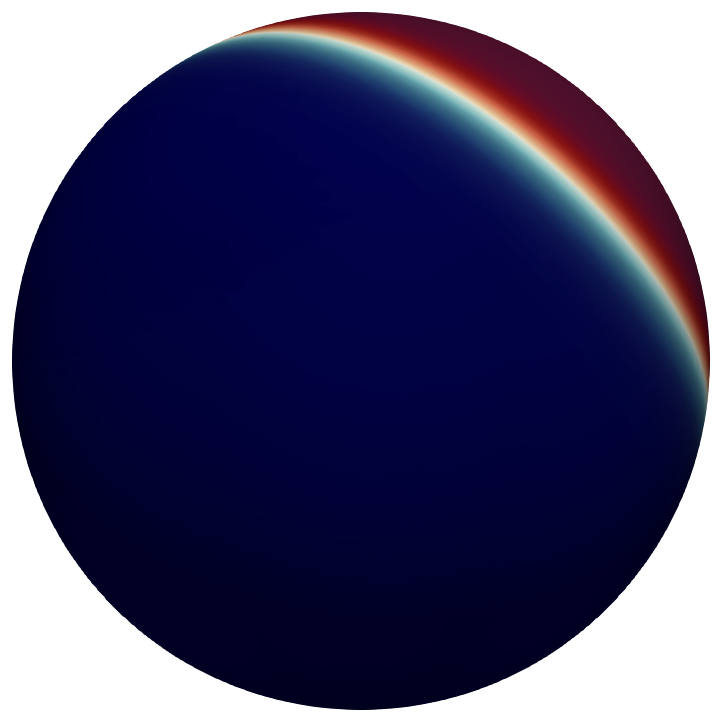}
			\put(30,102){\small{$t = 10$}}
		\end{overpic}
		\begin{overpic}[width=.13\textwidth,grid=false]{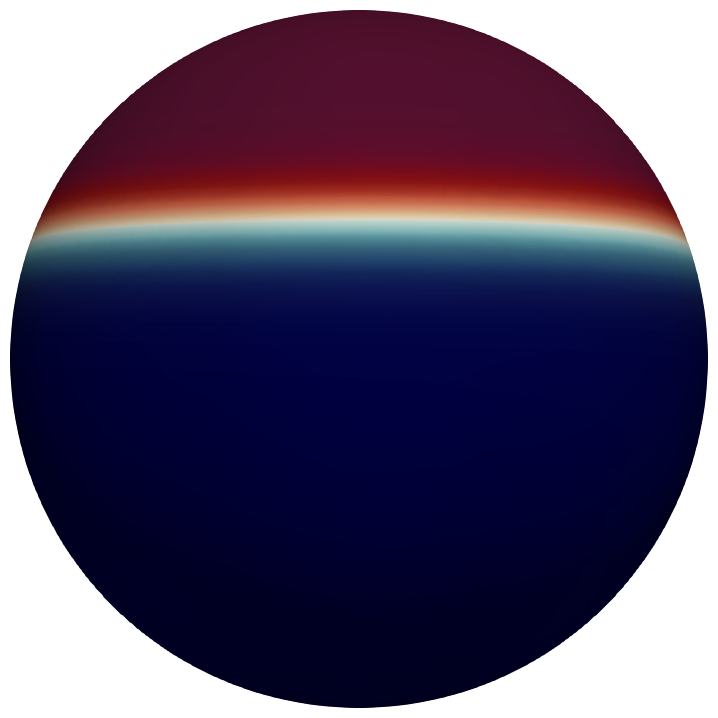}
			\put(30,102){\small{$t = 15$}}
		\end{overpic}\\
		
		\vskip 12pt
        \begin{overpic}[width=.13\textwidth,grid=false]{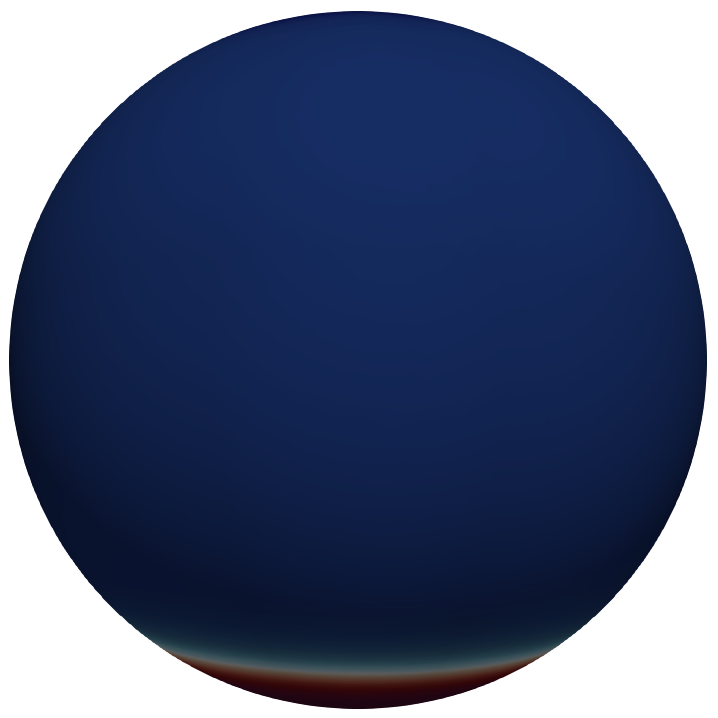}
			\put(30,102){\small{$t = 0$}}
			\put(-85,55){\small{PAT1}}
			\put(-108,37){\small{$a_D = 10.8$\%}}
		\end{overpic}
		\begin{overpic}[width=.13\textwidth,grid=false]{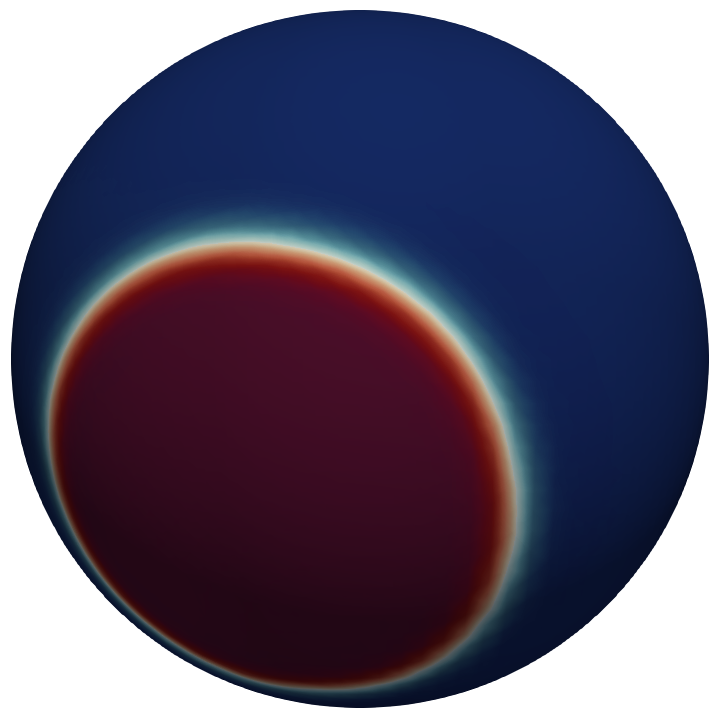}
			\put(30,100){\small{$t = 20$}}
		\end{overpic}
		\begin{overpic}[width=.13\textwidth,grid=false]{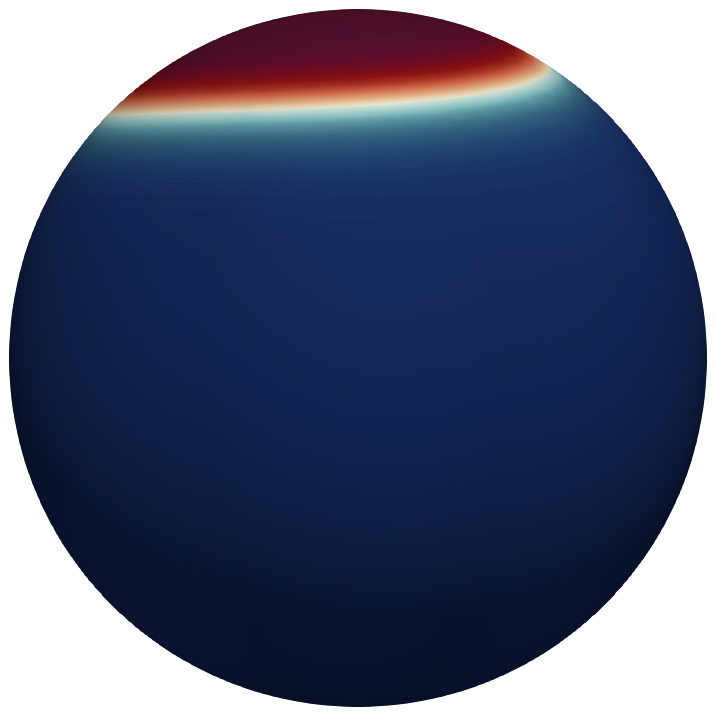}
			\put(30,102){\small{$t = 60$}}
		\end{overpic} \\
		\begin{overpic}[width=.5\textwidth,grid=false]{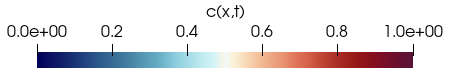}
		\end{overpic}

	\end{center}
	\caption{Snapshots from a simulation with the phase-separated PAT3 SUV (top), PAT2 SUV (center), and PAT1 SUV (bottom)
	at different times (min). 
	Red corresponds to the $L_o$ phase and blue to the $L_d$ phase. For each composition, the 
	$L_d$ phase is initially placed at the top of the SUV (first column).  The GUV, not seen in the figure, 
        is represented as a horizontal plane below the SUV.}
	\label{fig:numerical_snapshots}
\end{figure*}

Initially, the $L_d$ phase, which is the
phase with the majority of the positive charge, is placed opposite to the GUV, i.e., at the top of the SUV. 
See the first column in Fig.~\ref{fig:numerical_snapshots}. 
This is the worst-case scenario for fusion, 
as it will take the longest to get the $L_d$ phase to face 
the GUV and hence to initiate fusion. Recall that the majority of
the fusogenic lipids is in the $L_d$ phase (blue in Fig.~\ref{fig:numerical_snapshots}).
Fig.~\ref{fig:numerical_snapshots} shows simulation 
snapshots for the three compositions. We see that each SUV takes a different amount of time to have the 
$L_d$ phase face the model membrane. \revt{On} average (over 5 simulations), 
it takes 60 minutes PAT1 SUVs, 15 minutes for PAT2 SUVs and
6 minutes for PAT3 SUVs.
We take these times as a proxy for the promotion of fusion since they 
are the times needed to have the SUV in the optimal configuration for fusion, i.e., with the majority of the fusogenic lipids facing the GUV.
In average, a PAT1 SUV takes ten times longer 
than a PAT3 SUV to reorient its $L_d$ phase. Recall that the data used for Fig.~\ref{fig:fusion_frac} were acquired after
10 min of incubation. In that amount of time, the simulations predict that all PAT3 SUVs were in the optimal configuration 
for fusion, regardless of the initial position of the $L_d$ phase with respect to the GUV. In contrast, the PAT1 and PAT2 SUVs exposed to a GUV
in the worst-case scenario (i.e., $L_d$ phase opposite to the GUV) did not have sufficient time to have the $L_d$ phase face the GUV. 
Hence, the results in Fig.~\ref{fig:numerical_snapshots} provide
an explanation for the experimental data in Fig.~\ref{fig:fusion_frac}.

\section{\revo{Future perspectives}}\label{sec:fp}

\revo{
Preliminary experiments reported in \cite{wang2024} suggest a non-trivial interplay between 
charge density (i.e., the amount of DOTAP) and phase separation: 
charge density above a certain threshold inhibits phase separation. Thus, with the goal of assisting the design of phase-separated liposomes with fusogenic lipids in mind, one needs to extend the NSCH model to account for electrostatic charge and its interplay with phase separation and surface flow.}

\revo{
Electrostatic charges can lower miscibility temperatures in lipid membranes and alter the length scale of the phase-separated pattern~\cite{C4SM01089B,VequiSuplicy2010}. Thus, it is reasonable to hypothesize that the effects of electrostatic charges are long-range. A possible way to extend the NSCH model is to add a nonlocal term to the Ginzburg–Landau energy functional~\eqref{eq:f0}.
Heuristically, this term arises from considering long-range effects when deriving the system’s free energy using mean-field theory and statistical mechanics~\cite{Ohta1986,Oono1990,10.1063/1.5037727,bates1999}. The nonlocal term could be fractional Laplacian, resulting in the Otha-Kowasaki energy 
functional \cite{NISHIURA199531}. This is motivated by analyses in 1D \cite{NISHIURA199531, johnson2013branch, lessard2017rigorous} showing that this fractional term alters the threshold of phase separation,
increases the number of local minima, and yields metastable states with length-scales that depend on the strength of the nonlocal interaction. Obviously, suitable numerical algorithms would have to be designed for the resulting extension of the NSCH model. 
}

\revo{Moreover, future work for the design of phase-separated liposomes needs to account for the trigger of phase separation, since phase separation in liposomes should occur only close to the target cells, i.e., at the cancerous site. One strategy to control the onset
of phase separation consists in the use of polyethylene glycol (PEG) polymers linked to the liposome surface,
which inhibit phase separation by exerting pressure over the lipids. The acidic pH in the proximity of cancer cells cuts the links resulting in the separation of phases. To investigate this strategy, 
steric pressure exerted by PEG would have to be added to the NSCH model.}

\section{Conclusions}\label{sec:concl}

Phase-separated lipid vesicles serve as a powerful model system for studying the complex and dynamic behavior of biological membranes, particularly the spatial heterogeneity of physical properties across different membrane regions. Over the past two decades, these systems have attracted growing interest in biology, medicine, and nanotechnology.

This paper presented an overview of the full research cycle involved in modeling, simulating, and validating membrane dynamics. Beginning with the development of continuum-scale models, we described the surface Navier–Stokes equations for lipid flows, the surface Cahn–Hilliard equations for lateral phase separation, and their coupling in the surface NSCH system. We also included  a mechanical model accounting for membrane bending elasticity.
To connect theory with practice, we incorporated experimentally derived physical parameters into the models. 

We addressed the challenges of solving PDEs posed on surfaces using unfitted finite element methods. In particular, we presented TraceFEM, which extends surface PDEs into a volumetric neighborhood using a background mesh independent of the surface geometry, and  exemplified its application to the 
NSCH model.

Finally, we completed the cycle by validating the numerical predictions against wet-lab experiments for membranes of varying compositions. This direct comparison confirmed the NSCH model's ability to reliably capture phase behavior, even in the presence of  DOTAP --- a cationic lipid known to enhance membrane fusogenicity.

By including model development, parameter acquisition, numerical simulation, and experimental validation, this work reviewed a holistic approach to studying phase behavior in lipid membranes.

\section*{Acknowledgements}
All the experimental data used for validation were provided
by Dr.~Sheereen Majd and Dr.~Yifei Wang. Their work and 
availability for discussion were instrumental in the authors'
understanding of the complexity of lipid vesicles. 
The collaborative work with Drs.~Majd and Wang
was partially supported by US National Science Foundation (NSF) through grant DMS-1953535.
M.O. is also partially supported by National Science Foundation under Grants DMS-2309197 and DMS-2408978.

\section*{Author Declarations}
The authors have no conflicts to disclose.

\section*{Data Availability Statement}
The data that support the findings of this study are available from the corresponding author, A.~Q., upon reasonable request.

\nocite{*}
\bibliography{ref}

\end{document}